\documentclass{article}

\usepackage{amsmath, amsthm, amscd, amsfonts, amssymb, graphicx, color}
\usepackage[bookmarksnumbered, colorlinks, plainpages]{hyperref}
\usepackage{pspicture,pst-all,multido}
\theoremstyle{definition}

\theoremstyle{remark}

\numberwithin{equation}{section}

\newcommand{\SMALL}{\scriptsize}

\newcommand{\C}{\mathbb{C}}
\newcommand{\R}{\mathbb{R}}
\newcommand{\N}{\mathbb{N}}

\newcommand{\norm}[1]{\parallel\!\!#1\!\!\parallel}
\renewcommand{\d}{\operatorname{d}}
\newcommand{\e}{\operatorname{e}}
\renewcommand{\i}{\operatorname{i}}

\newcounter{envcount}%

\newenvironment{Def}%
{\vspace{\bigskipamount}\refstepcounter{envcount}\textbf{(\theenvcount)\enspace Definition.}}%
  {\vspace{\bigskipamount}}

{\vspace{\bigskipamount}\refstepcounter{envcount}\textbf{(\theenvcount)\enspace Summary.}}%
  {\vspace{\bigskipamount}}

{\vspace{\bigskipamount}\refstepcounter{envcount}\textbf{(\theenvcount)\enspace Example.}}%
  {\vspace{\bigskipamount}}

{\vspace{\bigskipamount}\refstepcounter{envcount}\textbf{(\theenvcount)\enspace Example}}%
  {\vspace{\bigskipamount}}

{\vspace{\bigskipamount}\refstepcounter{envcount}\textbf{(\theenvcount)\enspace Construction}}%
  {\vspace{\bigskipamount}}

{\vspace{\bigskipamount}\refstepcounter{envcount}\textbf{(\theenvcount)\enspace Theorem of Plancherel and P\'olya. }}%
  {\vspace{\bigskipamount}}

{\vspace{\bigskipamount}\enspace \textbf{Definition.}}%
  {\vspace{\bigskipamount}}

\newenvironment{The}%
{\vspace{\bigskipamount}\refstepcounter{envcount}\textbf{(\theenvcount)\enspace Theorem.}\itshape}%
  {\vspace{\bigskipamount}\upshape}

{\vspace{\bigskipamount}\refstepcounter{envcount}\textbf{(\theenvcount)\enspace Theorem}\itshape}%
  {\vspace{\bigskipamount}\upshape}
  
\newenvironment{Pro}%
{\vspace{\bigskipamount}\refstepcounter{envcount}\textbf{(\theenvcount)\enspace Proposition.}\itshape}%
  {\vspace{\bigskipamount}\upshape}

\newenvironment{Cor}%
{\vspace{\bigskipamount}\refstepcounter{envcount}\textbf{(\theenvcount)\enspace Corollary.}\itshape}%
  {\vspace{\bigskipamount}\upshape}

\newenvironment{Lem}%
{\vspace{\bigskipamount}\refstepcounter{envcount}\textbf{(\theenvcount)\enspace Lemma.}\itshape}%
  {\vspace{\bigskipamount}\upshape}

{\vspace{\bigskipamount}\refstepcounter{envcount}\textbf{(\theenvcount)\enspace Lemma and Example.}\itshape}%
  {\vspace{\bigskipamount}\upshape}

{\vspace{\bigskipamount}\refstepcounter{envcount}\textbf{(\theenvcount)\enspace Example.}\itshape}%
  {\vspace{\bigskipamount}\upshape}

{\vspace{\bigskipamount}\refstepcounter{envcount}\textbf{(\theenvcount)\enspace Example}\itshape}%
  {\vspace{\bigskipamount}\upshape}

  {\vspace{\bigskipamount}}

  {\vspace{\bigskipamount}}

{\vspace{\bigskipamount}\refstepcounter{envcount}\textbf{(\theenvcount)\enspace Interpretation.}\itshape}%
  {\vspace{\bigskipamount}\upshape}

{\vspace{\bigskipamount}\refstepcounter{envcount}\textbf{(\theenvcount)\enspace}}%
  {\vspace{\bigskipamount}}

\theoremstyle{definition}
\swapnumbers

\begin{document}
\setcounter{page}{1}
\pagenumbering{arabic}

\begin{center}
{\Large Spectral Representations of the Wiener-Hopf Operator for the sinc Kernel and some Related  Operators}\\  

\vspace{1cm}
Domenico P.L. Castrigiano\\
Technischen Universit\"at M\"unchen, Fakult\"at f\"ur Mathematik, M\"unchen, Germany\\

\smallskip

{\it E-mail address}: {\tt
castrig\,\textrm{@}\,ma.tum.de
}
\end{center}

\begin{quote} The spectral representation  of the Wiener-Hopf  operator $K$ with kernel $\frac{1}{\pi}\operatorname{sinc}$ is given determining explicitly the Hilbert space isomorphism, which transforms $K$ into the multiplication operator by the identity on $L^2(0,1)$. Several related integral operators are studied. A close relationship of $K$ to the finite Hilbert transformation is revealed  yielding  the spectral representation of the latter.  This is of particular interest as it concerns a general  feature of self-adjoint Wiener-Hopf operators  \cite{C20}.\\

\end{quote}

\section{Introduction} 
We study the integral operator $K$ on $L^2(\R_+)$, $\R_+:=\{x\in\R: x>0\}$ given by 
\begin{equation}\label{WHO}
(Kg)(x)=\int_0^\infty \textrm{\SMALL{$\frac{1}{\pi}$}}\operatorname{sinc}(x-y)\,g(y)\d y
\end{equation}

deriving its spectral representation. The sinc function or cardinal sine function $\operatorname{sinc}(x)=\frac{\sin(x)}{x}$ for $x\in\R\setminus \{0\}$, $\operatorname{sinc}(0)=0$   (also called sampling function or slit function)   frequently arises in signal processing and the theory of Fourier transforms. It equals $\operatorname{j}_0$, which is
  the zeroth order spherical Bessel function of the first kind.  $K$ is  an operator of Wiener-Hopf type. It is not  singular  as $\operatorname{sinc}$ is continuous at $0$.
A well-known method  to achieve the  concrete  spectral representation of the operator requires the solution of the corresponding homogeneous integral equation (\ref{HWHE}).
The  original Wiener-Hopf method \cite{WHP31} does not apply   as   $\frac{1}{\pi}\operatorname{sinc}(x)\exp(\lambda|x|)$, $x\in\R$ is not square integrable for any $\lambda>0$.  
  Since the kernel $\frac{1}{\pi}\operatorname{sinc}$ is not even integrable and the symbol $1_{[-1,1]}=\sqrt{2\pi}\mathcal{F}^{-1}\big(\frac{1}{\pi}\operatorname{sinc}\big)$ is not continuous,
  more general theories as Krein \cite{K62}, Talenti \cite[Teorema 3.1.2, 3.1.3]{T73} do not apply, either. However, as we will see, one may  successfully  introduce the factor $\exp(-\varepsilon|\cdot|)$, $\varepsilon>0$ to ensure the exponential decay of the kernel, and study the limit $\varepsilon\to 0$.
An alternative approach to the spectral representation of $K$ applies the theory by Rosenblum 
 \cite{R61} by turning to an isomorphic Toeplitz operator, see sec.\,\ref{TOIK} for details.\\
\hspace*{6mm}
We came across $K$  when  analyzing   the causal localization  for a Weyl fermion. In sec.\,\ref{CLOWF} it is shown that the localization operator for the unit ball is Hilbert space isomorphic 
to the countably infinite orthogonal sum of copies of $K$. In the course of the proof some useful formulae concerning the Bessel functions of the first kind and the related Hankel transformations are deduced. --- Exploiting the Hilbert space isomorphisms transforming Wiener-Hopf operators into Toeplitz operators (\ref{MTRTOWHO}) we obtain in sec.\,\ref{TOIK} the spectral representation of the Toeplitz operators for the indicator functions of the arcs on the unit circle, already derived in \cite{R61}.  The integral transformations effecting the spectral representation are studied in (\ref{GOLBV}), (\ref{SGOLBV}). They map bases of known orthonormal functions onto each other. --- Particularly interesting is the close relationship of $K$ to the finite  
Hilbert transformation revealed in sec.\,\ref{MEFHT}, which allows to deduce the spectral representation of the latter \cite{KP59}. In \cite{C20} it is shown that generally a natural isometry exists transforming a self-adjoint Wiener-Hopf operator into a singular integral operator of the well-studied class of Hilbert transformation type operators.    
 --- Wiener-Hopf  operators for symbols which are related by  change of variable due to a real fractional linear transformation 
are unitarily equivalent. Actually the covariance in (\ref{CVFLTWH}) holds. Using this, in (\ref{HTWHRK}) several Wiener-Hopf operators and singular integral operators are transformed isomorphicly    into $2K-I$.

\section{Spectral Representation of $K$}

Let us call a solution $u$ of the homogeneous Wiener-Hopf equation 
\begin{equation}\label{HWHE}
s\,u(x)- \int_0^\infty k(x-y)u(y)\d y=0, \quad x>0
\end{equation}
a generalized eigenfunction or  \textbf{spectral function} of the respective Wiener-Hopf operator for the generalized eigenvalue or  \textbf{spectral value}   $s$. 
For the kernel $\frac{1}{\pi}\operatorname{sinc}$ we will show that every $s\in ]0,1[$ is a spectral value and find the corresponding  spectral functions $q_s$. Moreover, these are shown to  give rise to the kernel $q:]0,1[\times \R_+\to \C$, $q(s,x):=n(s)q_s(x)$ for the normalizing factor $n(s)=\frac{\i}{\pi}(2s)^{-1/2}$ such that the integral operator 
$(Vh)(x):=\operatorname{l.i.m.}\int_{\downarrow 0}^{\uparrow1}q(s,x)h(s)\d s$ determines a Hilbert space isomorphism from $L^2(0,1)$ onto $L^2(\R_+)$ diagonalizing $K$, i.e.,
$$\big(V^{-1}K\,Vh\big)(s)=s\,h(s)$$
Hence the spectrum of $K$ is $[0,1]$, it is simple,  and  $K$ is absolutely continuous. 
As known, the latter  property have all bounded self-adjoint Wiener-Hopf operators.

\subsection{First properties of K}\label{FPOK}
 $\mathcal{F}$ denotes the Fourier transformation on  $L^2(\R)$. Recall  $\mathcal{F}f(x)=(2\pi)^{-1/2}\int_\R \e  ^{-\i  xy} f(y) \d y$ for $x\in\R$ if $f\in L^2\cap L^1$ and that in general $\mathcal{F}f$ is the limit $\operatorname{l.i.m.}_{A,A'\to\infty}\mathcal{F}(1_{[-A,A']}f)$  in $L^2$. 
For any $\kappa\in L^\infty(\R)$ let $M(\kappa)$ denote the multiplication operator by $\kappa$ on $L^2(\R)$, which is normal and bounded with norm $\norm{\kappa}_\infty$.  We use the projection $P_+:L^2(\R)\to L^2(\R_+)$, $P_+f(x):=f(x)$ and its adjoint, the injection $P_+^*:L^2(\R_+)\to L^2(\R)$, $P_+^*g(x)=g(x)$ if $x>0$ and $=0$ otherwise. Note $P_+^*P_+=M(1_{\R_+})$ and $P_+P_+^*=I_{L^2(\R_+)}$.

\begin{Def}
Let $\kappa\in L^\infty(\R)$. Then the bounded operator $$W_\kappa:=P_+\mathcal{F}M(\kappa)\mathcal{F}^{-1}P_+^*$$ on $L^2(\R_+)$ is called the Wiener-Hopf operator with symbol $\kappa$. Occasionally  we write $W(\kappa)$ instead of $W_\kappa$.
\end{Def} 

Obviously, $\norm{W_\kappa}\le \norm{\kappa}_\infty$, and $W_\kappa$ is self-adjoint for real $\kappa$ with the spectrum contained in $[\operatorname{ess-inf} \kappa, \operatorname{ess-sup}\kappa]$. Actually, as known,  the spectrum equals this interval.

 \begin{Lem}\label{WHMCO}
 Consider the two cases
 \begin{itemize}
 \item[\emph{(a)}] $\kappa\in L^\infty(\R)\cap L^2(\R)$ and $k:=(2\pi)^{-1/2}\mathcal{F}\kappa$
  \item[\emph{(b)}] $\kappa=\int\e  ^{\i  (\cdot) y}k(y)\d y$ for some $k\in L^1(\R)$
 \end{itemize}
 In either case   the  operator $C(k):=\mathcal{F}M(\kappa)\mathcal{F}^{-1}$ on $L^2(\R)$
 is bounded with $\norm{C(k)}=\norm{\kappa}_\infty$ and satisfies 
 $$\big(C(k)f\big)(x)=\int_{-\infty}^\infty k(x-y)f(y)\d y$$ 
 In case \emph{(a)} the integral exists  for every $x\in\R$ and determines a continuous function  bounded by $\norm{k}_2\norm{f}_2$. 
 \end{Lem}\\
 {\it Proof.} (a) Write $\mathcal{F}^{-1}f=\operatorname{l.i.m.}_{A\to\infty}F_A$ for $F_A(z):=(2\pi)^{-1/2}\int_{-A}^A \e  ^{\i zy} f(y) \d y$. Then  by continuity 
$C(k)f=\mathcal{F}M(\kappa)\mathcal{F}^{-1}f=\operatorname{l.i.m.}_{A\to\infty}\mathcal{F}M(\kappa)F_A$. Furthermore,
 $\mathcal{F}M(\kappa)F_A=\operatorname{l.i.m.}_{B\to\infty}G_{A,B}$ with $G_{A,B}(x):= (2\pi)^{-1/2}\int_{-B}^B \kappa(z)\e  ^{-\i xz}\big((2\pi)^{-1/2}\int_{-A}^A \e  ^{\i zy} f(y) \d y\big)\d z$ for $x\in\R$. Then by  Fubini's theorem, $G_{A,B}(x)=\int_{-A}^A K_{x,B}(y)f(y)dy$ for $K_{x,B}(y):=\frac{1}{2\pi}\int_{-B}^B\kappa(z)\e  ^{-\i  z(x-y)}\d z$. As $\operatorname{l.i.m.}_{B\to\infty}K_{x,B}=k\big(x-(\cdot)\big)$, one has $\lim_{B\to\infty}G_{A,B}(x)=\int_{-A}^A k(x-y)f(y)\d y$ for every $x$. Now note $\int|k(x-y)f(y)|\d y\le \norm{k}_2\norm{f}_2$ by Cauchy-Schwarz inequality, whence $\big(C(k)f\big)(x)=\int k(x-y)f(y)\d y$. Moreover, $\phi(x):=\int k(x-y)f(y)\d y$  exists for every $x$ and $2\pi|\phi(x+h)-\phi(x)|^2\le\norm{f}_2^2\int|\big(\e  ^{-\i  xh}-1\big)\kappa(x)|^2\d x\to 0$ for $h\to 0$ by dominated convergence. \\
\hspace*{6mm}  (b) In view of    \cite[(21.31),\,(21.32)]{HS69}   it suffices  to show $\kappa\,\mathcal{F}^{-1}f=\mathcal{F}^{-1}(k\star f)$ for $f\in L^2\cap L^1$. Now,  $k\star f\in L^2\cap L^1$ by \cite[(21.31),\,(21.32)]{HS69}
and hence $\mathcal{F}^{-1}(k\star f)=\kappa\, \mathcal{F}^{-1}f$ by \cite[(21.41)]{HS69}.\qed

In the cases treated in  (\ref{WHMCO}) the Wiener-Hopf operator $W_\kappa$  is the trace  on the subspace $L^2(\R_+)$ of the convolution operator $C(k)$ on $L^2(\R)$ and hence has kernel $k$. We apply now (\ref{WHMCO})(a).

\begin{Cor}\label{FPK} As 
$\frac{1}{\pi}\operatorname{sinc}=(2\pi)^{-1/2}\mathcal{F}1_{[-1,1]}$ one has
$$K=W(1_{[-1,1]})$$  
with its spectrum  contained in $[0,1]$. Note $\norm{\frac{1}{\pi}\operatorname{sinc}}_2=\pi^{-1/2}$.
For $g\in L^2(\R_+)$ the right hand side of \emph{(\ref{WHO})} exists for every $x\in\R$ and $Kg\in L^2(\R_+)$ is bounded continuous  with $\norm{Kg}_2\le\norm{g}_2$ and $\norm{Kg}_\infty\le\pi^{-1/2}\norm{g}_2$.
\end{Cor}

\subsection{An auxiliary Wiener-Hopf operator}\label{AAWHO} 
For $\varepsilon>0$   introduce the Wiener-Hopf operator $W_{\kappa_\varepsilon}$ with  kernel $$k_\varepsilon=\textrm{\SMALL{$\frac{1}{\pi}$}}\e  ^{-\varepsilon |\cdot|}\operatorname{sinc}$$ As $k_\varepsilon$ vanishes exponentially at infinity the Wiener-Hopf method \cite[11.17]{T48} applies and yields a spectral function  of $W_{\kappa_\epsilon}$ for every  value $s\in]0,1[$.  Then the limit $\varepsilon\to 0$ results in a family $q_s$, $s\in]0,1[$ of functions,  which in (\ref{SFKP}) are shown  to be  spectral functions  of $K$.

Following \cite[11.17]{T48} consider the symbol $\kappa_\varepsilon(w)= \int_{-\infty}^\infty\e  ^{\i  wy} k_\varepsilon(y) \d y$ for complex $w:=u+\i  v$, $|v|<\varepsilon$: 
\begin{eqnarray*}
2\pi\operatorname{Re} \kappa_\varepsilon(w)&=&\arctan\textrm{\SMALL{$\left(\frac{u+1}{\varepsilon-v}\right)$}}-\arctan\textrm{\SMALL{$\left(\frac{u-1}{\varepsilon-v}\right)$}}+\arctan\textrm{\SMALL{$\left(\frac{u+1}{\varepsilon+v}\right)$}}-\arctan\textrm{\SMALL{$\left(\frac{u-1}{\varepsilon+v}\right)$}}\\
2\pi\operatorname{Im} \kappa_\varepsilon(w)&=&\textrm{\SMALL{$\frac{1}{2}$}}\ln\textrm{\SMALL{$\left(\frac{(u-1)^2+(\varepsilon-v)^2}{(u+1)^2+(\varepsilon-v)^2} \cdot\frac{(u+1)^2+(\varepsilon+v)^2}{(u-1)^2+(\varepsilon+v)^2}\right)$}}
\end{eqnarray*}
In view of the spectrum of $K$ (see (\ref{FPK})) restrict at once $s\in\,]0,1[$. For $\varepsilon>0$ sufficiently small, i.e., $\varepsilon<\cot(\frac{\pi}{2}s)$, the function to be factorized $1-\frac{1}{s}\kappa_\varepsilon$ has  exactly  two  zeros, namely at $\pm\, u_{\varepsilon,s}$ for $u_{\varepsilon,s}:=\big(1-\varepsilon^2+2\varepsilon\cot(\pi s)\big)^{1/2}$, which are real and simple. Put the cut of $\ln$ along the negative real axis. By the factor $R_ {\varepsilon,s}(w):=\frac{w^2+\varepsilon^2}{w^2-u_ {\varepsilon,s}^2}$ the product $S_{\varepsilon,s}:=(1-\frac{1}{s}\kappa_\varepsilon)R_ {\varepsilon,s}$  has no zeros in the strip $\{|v|<\varepsilon\}$ and the variation of $\ln\circ \,S_{\varepsilon,s}$ along the strip is zero. Let
\begin{equation}\label{AWHOF}
T_{\varepsilon,s}(w):=\textrm{\SMALL{$\frac{1}{2\pi\i  }$}}\int_{-\infty-\i  b}^{\infty-\i  b}\frac{\ln S_{\varepsilon,s}(z)}{z-w}\d z,\quad Q _{\varepsilon,s}(w):=\e  ^{-T_{\varepsilon,s}(w)}\frac{w+\i  \varepsilon}{w^2-u_{\varepsilon,s}^2}
\end{equation}
with $0<b<\varepsilon$ for $v>-b$. Then
\begin{equation}
q _{\varepsilon,s}(x):=(2\pi)^{-1/2}\int_{-\infty+\i  b}^{\infty+\i  b}
\e  ^{-\i  xw}Q _{\varepsilon,s}(w)\d w
\end{equation}
is a spectral function of $W_{\kappa_\epsilon}$ for the spectral value $s\in]0,1[$ (and $q _{\varepsilon,s}(x)=0$ for $x<0$). Moreover,  within the set of  locally integrable functions being $\mathcal{O}(\e  ^{cx})$ for some $c<\varepsilon$ as $x\to \infty$,
$q _{\varepsilon,s}$ is unique up to a complex constant (depending on $s$). Indeed, in case of even kernel a possible polynomial, by which $Q _{\varepsilon,s}$ has to be multiplied, has degree not exceeding $\frac{n}{2}-1$ with $n$ the number of zeros, see \cite[Theorem XVI in Paley,\,Wiener]{WHP31}. 

\subsubsection{The limit $\varepsilon\to 0$} 
Let $w=u+\i  v\in\C$.
From (\ref{AWHOF}) one obtains $\lim_{\varepsilon\to 0}T_{\varepsilon,s}(w)=\frac{1}{2}\ln\big(\frac{w^2}{w^2-1}\big)-\frac{\i  }{2\pi}\ln(\frac{1}{s}-1)\ln\big(\frac{w-1}{w+1}\big)$ and accordingly we put 
\begin{equation}\label{QS}
Q_s(w):=
\textrm{\SMALL{$\frac{w}{w^2-1}$}}
\exp\textrm{\SMALL{$\left(-\frac{1}{2}\ln\big(\frac{w^2}{w^2-1}\big)\right)$}}
\exp\textrm{\SMALL{$\left(\frac{\i  }{2\pi}\ln\big(\frac{1}{s}-1\big)\ln\big(\frac{w-1}{w+1}\big)\right)$}}
\end{equation}

One easily verifies  that $Q_s$ is holomorphic on $\C\setminus[-1,1]$ without zeros, has a cut along $[-1,1]$ and a singularity at $1$ and $-1$, and that
\begin{equation}\label{EQSW}
|Q_s(w)|\le C_s|w^2-1|^{-1/2} \le C_s|u^2-1|^{-1/2} \textrm{\, for \,} C_s:=\big(\textrm{\SMALL{$\frac{1}{s}-1$}}\big)^{1/2}
\end{equation}
Also one checks that $Q_s^\pm(u):=\lim_{v\to\pm 0}Q_s(u+\i  v)$ exists with
\begin{equation}\label{LQS}
Q_s^\pm(u)=m_s(u)^{\mp 1}|u^2-1|^{-1/2}\exp\textrm{\SMALL{$\left(\frac{\i  }{2\pi}\ln\big(\frac{1}{s}-1\big)\ln \big|\frac{u-1}{u+1}\big|\right)$}}
\end{equation}
where $m_s(u):=\operatorname{sgn}(u)$ if $|u|>1$ and $m_s(u):=\i  \big(\frac{1}{s}-1\big)^{1/2}$ if $|u|<1$. By (\ref{EQSW}) one has
\begin{equation}\label{LQSI}
Q_s^\pm\in L^{p}(\R) \textrm{ \,for \,} 1<p<2 \textrm{\; \;and \;\;} 1_{[-1,1]}Q_s^\pm\in L^{p}(\R)\textrm{ \,for \,} 1\le p<2
\end{equation}

Finally, for $b>0$ and $x\in\R\setminus\{0\}$ define $q_s^\pm(x)\in \C$ by the Cauchy integral
\begin{equation}\label{QSP}
q _s^\pm(x):=(2\pi)^{-1/2}\int_{-\infty\pm\i  b}^{\infty\pm\i  b}
\e  ^{-\i  xw}Q _s(w)\d w
\end{equation}

\subsubsection{Remark on the factorization} Note that (\ref{LQS}) yields    
 the factorization 
\begin{equation}
1-\textrm{\SMALL{$\frac{1}{s}$}} 1_{[-1,1]}(u)=Q^-_s(u)/Q^+_s(u), \quad u\in\R\setminus\{-1,1\}
\end{equation}
 for every $s\in]0,1[$, where by (\ref{QS}) $Q^-_s$ and $1/Q^+_s$ is the limit of a  holomorphic function  on the lower and upper half-plane, respectively. Such a factorization of $1-\frac{1}{s}\kappa$ is at the core of  the method for solving (\ref{HWHE}) by Wiener and Hopf and generalized by Krein. In general it is available,   if e.g. $\kappa$ is holomorphic in a strip around the real axis and uniformly square integrable on lines parallel to the real axis (see e.g. \cite[11.17 Lemma]{T48}), or more generally if $\kappa$ is $\alpha$-H\"older continuous ($0<\alpha<1$) \cite[Teorema 2.1.2]{T73}, or if $\kappa$ is the Fourier transform of an integrable function \cite{K62}, see e.g. \cite[Teorema 2.1.4]{T73}.

\subsection{Spectral functions of $K$} We are going to show that $q_s:=q_s^+|_{\R_+}$  are spectral  functions of $K$. We start from their definition  (\ref{QSP}) and forget about their origin in  sec.\,\ref{AAWHO}. For the proof and the ensuing  results it is convenient to consider  also the functions $q_s^-$. Recall $s\in\,]0,1[$ and that $\ln$ denotes the principal branch of the logarithm. 

\begin{Pro}\label{MPQS} The integral defining $q_s^\pm(x)$ in \emph{(\ref{QSP})} exists 
 in a uniform manner for $s\in[\delta,1[$\,and $x\in]-\infty,-\delta]\cup[\delta,1/\delta]$ for every $\delta\in]0,1[$, and satisfies
\begin{itemize}
\item[\emph{(i)}] $q_s^\pm(x)=0$ for $\pm\, x<0$
\item[\emph{(ii)}] $q_s^+(x)=(2\pi)^{-1/2}\,\frac{1}{s}\int_{-1}^1\e  ^{-\i  xu}Q_s^+(u)\d u$  for $x>0$,\\
$q_s^-(x)=(2\pi)^{-1/2}\,\frac{1}{1-s}\int_{-1}^1\e  ^{-\i  xu}Q_s^-(u)\d u$  for $x<0$
\item[\emph{(iii)}] $q_s^\pm(x)=(2\pi)^{-1/2}\int_{ -\infty}^{ \infty}\e  ^{-\i  xu}Q_s^\pm(u)\d u$  for $x\ne 0$,  more precisely\\ $(2\pi)^{-1/2}\int_{ -A}^A\e  ^{-\i  xu}Q_s^\pm(u)\d u\to q_s^\pm(x)$ as $A\to \infty$ uniformly with respect to $x$ and $s$  in any compact subset of $\R\setminus\{0\}$ and $]0,1[$, respectively
\item[\emph{(iv)}] $\big|q_s^\pm(x)\big|\le(\pi/2)^{1/2}\big(s(1-s)\big)^{-1/2}$ for $x\ne 0$, and
 $q_s^\pm\in L^p(\R)$ for $p>2$ 
\end{itemize}
\end{Pro}
{\it Proof.} We prove the assertions for $q_s^+$. The results for $q_s^-$ follow quite analogously. Put $I(x,s,A):=\int_{-A +\i  b}^{A+\i  b}\e  ^{-\i  xw}Q _s(w)\d w$ for $A>0$.\\

\hspace*{6mm}
(i) By (\ref{EQSW}) one gets $\big|\e  ^{-\i  xw}Q_s(w)\big|  \le C_s \e  ^{xv}(u^2+v^2-1)^{-1/2}$ for $x\in\R$. Now let $x<0$ and consider the contours $L_1,\dots, L_4$ shown in Fig.\,1.

\begin{pspicture}(-3.5,-1)(3.5,4)
\psline{->}(0,-1)(0,3.5)
\psline{->}(-3.1,0)(3.1,0)
\psframe[linewidth=0.5mm](-2.5,0.8)(2.5,2.5)
\put(-2.5,0.1){-A}
\put(2.5,0.1){A}
\put(0.1,2.6){iA}
\put(0.1,0.95){ib}
\put(1.5,0.5){$\rightarrow$}

\put(-1.5,0.95){$L_1$}
\put(2.6,1.4){$L_2$}
\put(-1.5,2.6){$L_3$}
\put(-3,1.4){$L_4$}
\put(-3.1,-0.7){Fig. 1}
\end{pspicture}

\vspace{2mm}
Then $\left|\int_{L_2}\e  ^{-\i  xw}Q_s(w)\d w\right|\le C_s \int_0^A \e  ^{xv}(A^2-1)^{-1/2}\d v=C_s(A^2-1)^{-1/2}\frac{1}{x}\big(\e  ^{xA}-1\big)\to 0$ as $A\to\infty$. Similarly, the integral along $L_4$  vanishes for $A\to \infty$. Also $\left|\int_{L_3}\e  ^{-\i  xw}Q_s(w)\d w\right|\le C_s\int_{-A}^A   \e  ^{xA}(u^2+A^2-1)^{-1/2}\d u\le C_s \e  ^{xA}\int_{-A}^A1\d u=2C_sA\e  ^{xA}\to 0$ as $A\to \infty$. Now let $L$ be the closed contour composed by $L_1,\dots, L_4$. Since $\e  ^{-\i  x(\cdot)}Q_s$ is holomorphic in the upper half-plane, one has $0=\int_L\e  ^{-\i  xw}Q_s(w)\d w\to 
(2\pi)^{1/2}q _s^+(x)$ as $A\to\infty$. ---  
The estimates show that the limit $I(x,s,A)\to 0$ for $A\to\infty$ 
occurs uniformly with respect to  $x$ and $s$ if $x\le -\delta$ and $\delta\le s<1$ for some $\delta>0$.\\

(ii) Let $x>0$. By holomorphy  the integrals of $\e  ^{-\i  x(\cdot)}Q_s$ along $L$ and $L'$ are equal for every $A>1$. Hence by the same kind of estimations as used in the proof of (i) one obtains 
 $I(x,s,A)
\to (2\pi)^{1/2}q _s^+(x)$ for $A\to\infty$ and thus $(2\pi)^{1/2}q _s^+(x)=\int_{L'}\e  ^{-\i  xw}Q_s(w)\d w$. --- Note that the limit occurs uniformly with respect to  $x$ and $s$, for $x$ in a compact set $\subset\R_+$ and $s\in[\delta,1[$ for $\delta>0$.

\begin{pspicture}(-5,-5)(5,2)
\psline{->}(0,-4.5)(0,2)
\psline{->}(-4.6,0)(4.6,0)
\psframe[linewidth=0.5mm](-3.4,-3.4)(3.4,1)
\put(2,1.1){$\rightarrow$}
\psarc[linewidth=0.5mm]
(-1.8,0){5mm}{20}{340}
\psarc[linewidth=0.5mm]
(1.8,0){5mm}{-160}{160}
\put(-1.8,0.05){-1}
\put(-1.8,-0.04){$\centerdot$}
\put(1.9,0.07){1}
\put(1.8,-0.04){$\centerdot$}
\psline[linewidth=0.5mm](-1.35,0.2)(1.35,0.2)
\put(-0.8,0.3){$\rightarrow$}
\psline[linewidth=0.5mm](-1.35,-0.2)(1.35,-0.2)
\put(-2.8,-0.8){$L'_{-,\varepsilon}$}
\put(2.3,-0.8){$L'_{+,\varepsilon}$}
\put(-0.8,-0.6){$L'$}
\put(0.1,1.1){ib}
\put(0.1,-3.22){-iA}
\put(-2.8,-3.22){$L$}
\put(-3.95,0.1){-A}
\put(3.52,0.1){A}
\put(-4.6,-4.1){Fig. 2}
\end{pspicture}

Since the  singularities at $\pm1$ are integrable, the integrals along $L_{\pm,\varepsilon}$ vanish if $\varepsilon\to 0$. Hence $(2\pi)^{1/2}q _s^+(x)=\int_{-1}^1\e  ^{-\i  xu}\big(Q^+_s(u)-Q^-_s(u)\big)\d u$. By (\ref{LQS}) this yields the result. \\

(iii) Let $A>1$. Reasoning as in the proof of (i) one has $(2\pi)^{1/2}q _s^+(x)=\int_{L}\e  ^{-\i  xw}Q_s(w)\d w$  

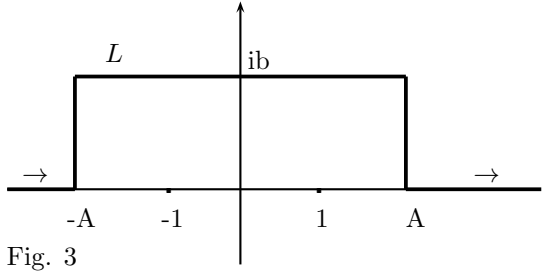
\begin{pspicture}(-3.5,-1)(3.5,3)
\psline{->}(0,-1)(0,2.5)
\psline(-3.1,0)(3.1,0)
\psline[linewidth=0.5mm](-2.2,0)(-2.2,1.5)
\psline[linewidth=0.5mm](2.2,0)(2.2,1.5)
\psline[linewidth=0.5mm](-2.2,1.5)(2.2,1.5)
\psline[linewidth=0.5mm](-3.1,0)(-2.2,0)
\psline[linewidth=0.5mm](2.2,0)(4,0)
\put(-1,-0.05){$\centerdot$}
\put(-1.05,-0.5){-1}
\put(1,-0.05){$\centerdot$}
\put(1,-0.5){1}
\put(-2.3,-0.5){-A}
\put(2.2,-0.5){A}
\put(0.1,1.6){ib}
\put(-1.8,1.7){$L$}
\put(-2.9,0.1){$\rightarrow$}
\put(3.1,0.1){$\rightarrow$}
\put(-3.1,-1){Fig. 3}
\end{pspicture}

\vspace{2mm}
(see Fig.\,3) independently of $b>0$. Note  
$I(x,s,A)
\to\int_{-A}^{A}\e  ^{-\i  xu}Q^+_s(u)\d u$ if $b\to 0$, since $C_s\e  ^{|x|}\big|u^2-1\big|^{-1/2}$ is integrable on $[-A,A]$ and  dominates $\e  ^{-\i  x(u+\i  b)}Q_s(u+\i  b)$  for $b\le 1$ by (\ref{EQSW}). Remember the uniformity of the limit of $I(x,s,A)$ for $A\to\infty$ in (i) and (ii). This implies the result. 

(iv) From (i) and (ii) and due to  (\ref{EQSW}) one gets $|q_s^+(x)| \le \frac{1}{s}(2\pi)^{-1/2}\int_{-1}^1|Q_s^+(u)|\d u=(\pi/2)^{1/2}\big(s(1-s)\big)^{-1/2}$.
 According to   \cite[Theorem 74]{T48} the results (iii)  and   (\ref{LQSI}) imply $q_s^+\in L^p(\R)$ for $p>2$.\qed

\begin{Cor} \label{EEQS}
$q_s^\pm|_{\R_\pm}$ has an entire extension.
\end{Cor}\\
{\it Proof.} Recall (\ref{MPQS})(ii) and (\ref{LQSI}). Hence by the Paley-Wiener theorem the Fourier-Laplace transform of $\frac{1}{s}1_{[-1,1]} Q_s^\pm$ is the asserted entire extension of $q_s^\pm|_{\R_\pm}$.\qed

\begin{Lem}\label{IKTQS}  Let $x\in\R$. The function $y\mapsto  \frac{1}{\pi}\operatorname{sinc}(x-y)q_s^\pm(y)$ on $\R$ belongs to $L^r(\R)$ for $r>\frac{2}{3}$ and hence  is integrable. The function $x\mapsto \frac{1}{\pi}\int_{-\infty}^{\infty} \operatorname{sinc}(x-y)q_s^\pm(y)\d y$ on $\R$ is continuous.
\end{Lem}\\
{\it Proof.} It suffices to consider the case $q_s^+$. Let $k$ equal $\operatorname{sinc}$ or more generally any $k\in\cap_{p>1} L^p(\R)$, as e.g. $k_0:\R\to\R$, $k_0(x):=1$ if $|x|\le 2$ and $k_0(x):=\frac{2}{|x|}$ if $|x|>2$.\\
\hspace*{6mm}
 By H\"older's inequality   $\int |q_s^+(x)|^r|k(x-y)|^r\d y\le \big(\int |q_s^+(x)|^{rp''}\d y\big)^{1/p''} \big(\int |k(x-y)|^{rq''}\d y\big)^{1/q''}$ with $1/p''+1/q''=1$. We choose $\frac{1}{r}<q''<\frac{3}{2}$ for $r>\frac{2}{3}$. Then $rq''>1$ and $1-r<1-\frac{1}{q''}=\frac{1}{p''}<\frac{1}{3}$. The latter implies $rp''>2$. Thus by (\ref{MPQS})(iv)  the right hand side of H\"older's inequality is finite. This implies the first part of the assertion.\\
\hspace*{6mm}
Let $f(x):=\int_{-\infty}^{\infty} \operatorname{sinc}(x-y)q_s^+(y)\d y$. Let $h\in\R$ with $|h|\le 1$. Check $\frac{1}{|t+h|}<\frac{2}{|t|}$ for $|t|>2$. Therefore $|\operatorname{sinc}(x+h-y)|<k_0(x-y)$
for all $x,y$. Hence $f(x+h)\to f(x)$ for $h\to 0$ by dominated convergence, whence the second part of the assertion.\qed

We turn to the main result of this section. Recall $q_s=q_s^+|_{\R_+}$ (see (\ref{QSP})). 

\begin{The}\label{SFKP}
Let $s\in\,]0,1[$. Then
$$\int_0^\infty \textrm{\SMALL{$\frac{1}{\pi}$}}\operatorname{sinc}(x-y)\,q_s(y)\d y\,=\,s\,q_s(x) \quad\forall \;x>0$$
Hence $q_s$ is a spectral function of $K$ for the spectral value $s$.
\end{The} \\
{\it Proof.} By (\ref{MPQS})(i) and (\ref{IKTQS}) the integral exists for every $x>0$.  By (\ref{MPQS})(iv) one has $q_s^+\in L^3(\R)$, $Q_s^+\in L^{3/2}(\R)$. Obviously, $k:=\frac{1}{\pi}\operatorname{sinc}\in L^{3/2}(\R)$. 
Because of (\ref{MPQS})(iii) the Fourier transform of $Q_s^+$ according to \cite[(4.1.1)]{T48} coincides a.e. with $q_s^+$. Hence the  Fourier transform of $q_s^+$ in the sense of \cite[(4.1.4]{T48}
coincides a.e. with $Q_s^+$. Recall $k(x)=(2\pi)^{-1/2}\int_{-\infty}^\infty  \e  ^{\i  xu}(2\pi)^{-1/2}1_{[-1,1]}(u)\d u$ for all $x\in\R$. Therefore the Fourier transform of $k$ in the sense of \cite[(4.1.1]{T48}
coincides a.e. with $(2\pi)^{-1/2}1_{[-1,1]}$.\\
\hspace*{6mm}
 Apply \cite[Theorem 76]{T48}. Thus $(2\pi)^{-1/2}\int_{-\infty}^\infty \e  ^{-\i  xu}(2\pi)^{-1/2}1_{[-1,1]}(u)Q_s^+(u)\d u=(2\pi)^{-1/2}\int_{-\infty}^\infty k(x-y)q_s^+(y)\d y$ a.e. Hence by (\ref{MPQS})(i),\,(ii) one has $s\,q_s(x)=\int_{0}^\infty k(x-y)q_s(y)\d y$ for almost all $x>0$. Actually equality holds for all $x>0$ since both sides are continuous functions of $x$ by (\ref{EEQS}) and (\ref{IKTQS}).\qed

\subsection{Hilbert space isomorphism associated with spectral functions of $K$}

As we will prove in (\ref{HSISFK}), the spectral functions $q_s$ of $K$ multiplied by the  normalization constants $n(s):=\frac{\i  }{\pi}(2s)^{-1/2}$ gives rise to a kernel $q$ for an integral operator, which determines an Hilbert space isomorphism $V$ from  $L^2(0,1)$ onto  $L^2(\R_+)$. 
The following definition considers both cases $q_s=q_s^+|_{\R_+}$ and $q_s^-|_{\R_-}$. 
\begin{equation}\label{KSFK}
q_\pm:]0,1[\times \R_\pm\to \C, \quad q_+(s,x):=n(s)q_s^+(x),\;\;q_-(s,x):=-n(1-s)q_s^-(x)
\end{equation}
By (\ref{MPQS})(ii) one has more explicitly 
\begin{equation}\label{EKSFK} 
q_+(s,x)=\textrm{\SMALL{$\frac{1}{2}$}}
\pi^{-3/2}s^{-1}(1-s)^{-1/2}\int_{-1}^1\e  ^{-\i  xu}
(1-u^2)^{-1/2}
\exp\textrm{\SMALL{$\left(\frac{\i  }{2\pi}\ln\big(\frac{1}{s}-1\big)\ln \big(\frac{1-u}{1+u}\big)\right)$}}\d u
\end{equation}
and $q_-(s,x)=q_+(1-s,-x)$. Note that $q_+$ and $q_-$ are real.

\begin{Lem}\label{SHSI} 
\begin{itemize} 
\item[\emph{(i)}] $R_\pm: L^2(0,1)\to L^2(\R,\C^2),\; (R_\pm h)(x):=(2\pi)^{1/2}\textrm{\SMALL{$\left(\begin{array}{c} l_\pm(s) \\ \pm\, l_\mp(s) \end{array}\right)$}}\,h(s)$ with $l_+(s):=s\,(1-s)^{1/2},\; l_-(s):=l_+(1-s)$ and $s:=\big(1+\e  ^{-2\pi x}\big)^{-1}$ for $x\in\R$, 
 is a linear isometry with range $\{f\in L^2(\R,\C^2):f_2(x)=\pm \e  ^{\mp\pi x}f_1(x) \operatorname{a.e.}\}.$ If $f$ is in the range of $R_\pm$, then  $(\R_\pm^{-1}f)(s)=(2\pi)^{-1/2}
 l_\pm(s)^{-1}\,f_1\big(\textrm{\SMALL{$-\frac{1}{2\pi}$}}\ln(\textrm{\SMALL{$\frac{1}{s}$}}-1)\big)$. \\

\item[\emph{(ii)}] $R:L^2(]0,1[,\C^2)\to L^2(\R,\C^2)$, $R\,h:=R_+h_1+R_-h_2$ is a Hilbert space isomorphism, $(R^{-1}f)(s)=(2\pi)^{-1/2}\textrm{\SMALL{$\left(\begin{array}{cc}(1-s)^{-1/2} & s^{-1/2} \\ s^{-1/2}  & -(1-s)^{-1/2} \end{array}\right)$}}f\big(\textrm{\SMALL{$-\frac{1}{2\pi}$}}\ln(\textrm{\SMALL{$\frac{1}{s}$}}-1)\big)$.\\

\item[\emph{(iii)}] $S: L^2(\R)\to L^2(\R,\C^2)$, $(Sf)(x):=\sqrt{2}\e  ^{x/2}\textrm{\SMALL{$\left(\begin{array}{c} (1+\e  ^x)^{-1}f\big(  (1-\e  ^x)  (1+\e  ^x)^{-1}\big)\\\i   \,(1-\e  ^x)^{-1}f\big(  (1+\e  ^x)  (1-\e  ^x)^{-1}\big)\end{array}\right)$}}$ is a Hilbert space isomorphism, $(S^{-1}f)(u)=\sqrt{2}\,|1-u^2|^{-1/2} \big(\i  \operatorname{sgn}(u)\big)^{1-j}f_j\big(\ln\big|\frac{1-u}{1+u}\big|\,\big)$ where $j=1$ if $|u|<1$ and $j=2$ if $|u|>1$.
\end{itemize}
\end{Lem}
{\it Proof.} The Hilbert space transformations all arise from change of variable and are easily checked.\qed
 
See sec.\,\ref{FPOK} for the projection $P_+$ and define $P_-$ analogously.

\begin{The}\label{HSISFK} $V_\pm:L^2(0,1)\to L^2(\R_\pm)$, $$V_\pm h:=\operatorname{l.i.m.}_{0<\alpha,\alpha'\to 0}\int_{\alpha}^{1-\alpha'}q_\pm(s,\cdot)h(s)\d s$$
 are Hilbert space isomorphisms with $$V_+^{-1}g=\operatorname{l.i.m.}_{A\to\infty}\int_0^Aq_+(\cdot,x)g(x)\d x, \quad V_-^{-1}g=\operatorname{l.i.m.}_{A\to\infty}\int_{-A}^0q_-(\cdot,x)g(x)\d x$$
 Moreover, the  formula 
$P_\pm^*\,V_\pm= \mathcal{F}^{-1}S^{-1} \big(\mathcal{F}^{-1}\oplus  \mathcal{F}^{-1}\big)R_\pm$ holds.\end{The}\\
{\it Proof.} Put $W_\pm:=\mathcal{F}^{-1}S^{-1} \big(\mathcal{F}^{-1}\oplus  \mathcal{F}^{-1}\big)R_\pm$ and  $f_{A,A'}:=(2\pi)^{-1/2}\int_{-A}^{A'}\e  ^{\i  (\cdot)y}\big(R_+h\big)(y)\d y$ for $h\in L^2(0,1)$ and $A,A'\ge 0$. Then $W_+h=\operatorname{l.i.m.}_{A,A'\to \infty}\mathcal{F}^{-1}S^{-1}f_{A,A'}$. \\
\hspace*{6mm}
For every $x\in\R$, one has 
$f_{A,A'}(x)=(2\pi)^{-1}\int_{\alpha}^{1-\alpha'}\textrm{\SMALL{$\e  ^{-\i  \frac{x}{2\pi}\ln\big(\frac{1}{s}-1\big)}$}}\textrm{\SMALL{$\left(\begin{array}{c} (1-s)^{-1/2} \\  s^{-1/2} \end{array}\right)$}}\,h(s)\d s$ by the change of variable $s:=\big(1+\e  ^{-2\pi y}\big)^{-1}$. Therefore, $G_{\alpha,\alpha'}(u):=\big(S^{-1}f_{A,A'}\big)(u)=\int_{\alpha}^{1-\alpha'}\overline{n(s)}\overline{Q_s^+(u)}\,h(s)\d s$ for every $u\in\R$. Put $J(x,s,B):=(2\pi)^{-1/2}\int_{-B}^B \e  ^{\i  xu}\overline{Q_s^+(u)}\d u$.\\
\hspace*{6mm}
Then $\mathcal{F}^{-1}G_{\alpha,\alpha'}=\operatorname{l.i.m.}_{B\to\infty}\int_{\alpha}^{1-\alpha'}\overline{n(s)}h(s)
J(\cdot,s,B)\d s$, 
where the order of integration has been inverted by Fubini's theorem since $\big|h(s) n(s)Q_s^+(u)\big|\le s^{-1}|h(s)||u^2-1|^{-1/2}$ (see (\ref{EQSW})) is integrable on $[\alpha,1-\alpha']\times[-B,B]$.\\
\hspace*{6mm}
It follows from (\ref{MPQS})(iii),(iv) that for $x\ne 0$ there is $B_x>0$ such that $|J(x,s,B)|\le |q_s^+(x)|+1\le (\pi/2)^{1/2}\big(s(1-s)\big)^{-1/2}+1$ for all $B\ge B_x$ and $s\in[\alpha,1-\alpha']$ and that $\lim_{B\to\infty}J(x,s,B)=\overline{q_s^+(x)}$. Thus 
$(\mathcal{F}^{-1}G_{\alpha,\alpha'})(x)=\int_{\alpha}^{1-\alpha'}\overline{n(s)}\overline{q_s^+(x)}h(s)\d s$ a.e. by dominated convergence. Recall (\ref{MPQS})(i) and that $q_+$ is real. Then this proves  $W_+=P_+^*\,V_+$. 
 In the same way one proves $W_-=P_-^*\,V_-$.\\
\hspace*{6mm}
By (\ref{SHSI})(ii), $W:L^2(]0,1[,\C^2)\to L^2(\R)$, $W\,h:=W_+h_1+W_-h_2$ is an isomorphism. This implies that $V_+$ and $V_-$ are isomorphisms, too.\\
\hspace*{6mm}
It remains to verify the integral representation for $V_+^{-1}$ and similarly that for $V_-^{-1}$. 
By the foregoing result on $W$ one has the formula $V_+^{-1}g=\big(W^{-1}P_+^*g\big)_1=R_+^{-1}\mathcal{F}^{(2)}S\mathcal{F}P_+^*g$ for $g\in L^2(\R_+)$, which is evaluated now. Let $A>0$ and $G_A:=(2\pi)^{-1/2}\int_0^A\e  ^{-\i  (\cdot)x}g(x)\d x$. Since $\int_0^A|g(x)|\d x<\infty$, $G_A$ is bounded on $\R$. Clearly $(SG_A)_1\in L^2(\R)$. Moreover, 
$(SG_A)_1(y)=\sqrt{2}\e  ^{y/2} (1+\e  ^y)^{-1}G_A\big(  (1-\e  ^y)  (1+\e  ^y)^{-1}\big)$, whence $\big|(SG_A)_1(y)\big|\le C\e  ^{y/2} (1+\e  ^y)^{-1}$ and hence $(SG_A)_1
\in L^1(\R)$. 
So $(2\pi)^{1/2}\big(\mathcal{F}(SG_A)_1\big)(x')=\int \e  ^{-\i  yx'}(SG_A)_1(y)\d y=\sqrt{2}\int_{-1}^1 
\exp\textrm{\SMALL{$\left(-\i  x'\ln \big(\frac{1-u}{1+u}\big)\right)$}}(1-u^2)^{-1/2}G_A(u)\d u$ by the new variable $u:=(1-\e  ^y)  (1+\e  ^y)^{-1}$. Now, Fubini's theorem obviously applies and, doing first the integral on $u$, one obtains at $x'=-(2\pi)^{-1}\ln\big(\frac{1}{s}-1\big)$ using (\ref{EKSFK}) the expression $\big(2\pi s^2(1-s)\big)^{-1/2}\big(\mathcal{F}(SG_A)_1\big)(x') =\int_0^Aq_+(s,x)g(x)\d x$, which equals  $\big(R_+^{-1}\mathcal{F}^{(2)}S\mathcal{F}G_A\big)(s)$
by  (\ref{SHSI})(i). Finally perform the limit in the mean for $A\to \infty$ on both sides.\qed

\hspace*{6mm}
See also (\ref{SGOLBV}) for a determining property of $V_+$.

\subsection{Spectral representation of $K$}\label{SRK} 

In the following put  $V:=V_+$.

\begin{The}
$\big(V^{-1}K\,V\,h\big)(s) = s \,h(s) \textrm{ \;a.e.\;}$ for every $h\in L^2(0,1)$.
\end{The}\\
{\it Proof.} Let $h\in L^2(0,1)$ vanish outside some interval $[\alpha,1-\alpha]\subset]0,1[$. The set of these functions is dense in   $L^2(0,1)$. By (\ref{FPK}),\,(\ref{HSISFK}), the proof of (\ref{HSISFK}), and (\ref{MPQS})(i) one has $KVh=P_+\mathcal{F}^{-1}M(1_{[-1,1]})\mathcal{F}P_+^*P_+\mathcal{F}^{-1}G_{\alpha,\alpha}=P_+\mathcal{F}^{-1}M(1_{[-1,1]})G_{\alpha,\alpha}$. Recall that $q_+$ is real. By  (\ref{EQSW}) Fubini's theorem applies so that by (\ref{MPQS})(ii)  $\big(KVh\big)(x)=\int_0^1q_+(s,x)\,s\,h(s)\d s$ for almost all $x>0$. Hence $\big(V^{-1}K\,V\,h\big)(s) = s \,h(s)$ a.e. by (\ref{HSISFK}). The result extends to all $h\in L^2(0,1)$ by continuity.\qed

So $K$ is Hilbert space isomorphic to the multiplication operator $M(\operatorname{id}_{[0,1]})$ by the identity  $\operatorname{id}_{[0,1]}$ on $L^2(0,1)$: 
\begin{equation}\label{KHSIJ}
 K\,=\,VM(\operatorname{id}_{[0,1]})V^{-1}
 \end{equation}

Recall  the spectral theorem (see e.g. \cite[12.23 Theorem]{R74}) for a self-adjoint operator $A$, by which  there exists a unique  spectral measure, i.e., a projection valued measure $E_A$ on the Borel sets of $\R$, satisfying $A=\int \lambda \d E_A(\lambda)$. Obviously, $E_{M(\operatorname{id}_{[0,1]})}(\Delta)h=1_{\Delta\cap [0,1]}h$, $h\in L^2(0,1)$, $\Delta\subset \R$ Borel set. So  $E_K(\Delta)=VE_{M(\operatorname{id}_{[0,1]})}(\Delta)V^{-1}$. 

\section{Some Operators related to $K$}

\subsection{Causal localization operator}\label{CLOWF} We came across $K$ (\ref{WHO})  when  analyzing    the causal localization  for a Weyl fermion, namely the localization operator $T^{\chi,\eta}(B)$ for the unit ball $B=\{x\in\R^3:|x|\le 1\}$, see (\ref{ALOWF}). One recalls that, as shown in \cite[sections 16, 23]{C18},  the two Weyl fermions ($\eta=+$) and their antiparticles ($\eta=-$) are the only massless relativistic quantum systems which satisfy the conditions imposed by causality. Here
$\chi\in\{+,-\}$ indicates the handedness of the particle, whence $m:=\chi\eta\frac{1}{2}$ its helicity.\\
\hspace*{6mm}
In the course of  the following considerations some useful formulae concerning  the Bessel functions of the first kind $\operatorname{J}_\nu$, $\nu\in\C$ and the related Hankel transformations $\mathcal{H}_\nu$, $\nu\in[-\frac{1}{2},\infty[$ are deduced, see (ii),\,(iii) and ($\alpha$),\,($\beta$) in (\ref{FAGKO}).

In position representation on $L^2(\R^3,\C^2)$ with $\R^3$ the coordinate space $T^{\chi,\eta}$ is the trace of the canonical projection valued measure  on the carrier space of the representation for the respective Weyl fermion, which in energy representation is simply $L^2(\R^2)\oplus\{0\}$ for $m=\frac{1}{2}$
and $\{0\}\oplus L^2(\R^2)$ for $m=-\frac{1}{2}$.\\
\hspace*{6mm}
 Let $T_m$ on $L^2(\R^3)$ denote $T^{\chi,\eta}(B)$ in the energy representation. For $\varphi\in L^2\cap L^1$  one has explicitly
\begin{equation}\label{CLOUB}
\big(T_m\varphi\big)(p)=(2\pi)^{-3/2}\int_{\R^3}D^{(1/2)}_{mm}(B_p^{-1}B_{p'})\big(\mathcal{F}1_B\big)(p-p')\varphi(p')\d p'
\end{equation}
Here $D^{(1/2)}$ is the identity representation of $SU(2)$ and $B_p\in SU(2)$ the helicity cross section satisfying $|p|B_pe_3=p$, $p\in\R^3$. Note $B_p=B_{p/|p|}$ for $p\ne 0$.     In the sequel identify $L^2(\R^3)$  by $p=s\,w$  with $L^2_\rho(\R_+)\otimes L^2_\omega(S^2)$, $\d \rho(s):=4\pi s^2\d s$, $\omega$ the normalized rotational invariant measure on the sphere $S^2$. Moreover consider the Hilbert space isomorphism $$\delta: L^2_\rho(\R_+)\to L^2(\R_+),\quad (\delta g)(s):=\sqrt{4\pi} s g(s)$$
 Recall that by the Peter-Weyl theorem the functions $w\to (2j+1)^{1/2}D^{(j)}_{m m'}(B^{-1}_w)$,  for $j\in |m|+\N_0$ and $m'\in\{-j,-j+1,\dots,j-1,j\}$, form an  ONB for $L^2_\omega(S^2)$. The orthogonal projection onto the subspace determined by fixed $j$ is given by 
$$(P^j_mf)(w)=(2j+1)\int_{S^2}D^{(j)}_{mm}(B^{-1}_wB_{w'})f(w')\d \omega(w')$$
Finally by \cite[7.14.2\,(32)]{E53} for $i\in\{0,1\}$ the functions $s\to (4n+2i+1)^{1/2}s^{-1/2}\operatorname{J}_{2n+i+\frac{1}{2}}(s)$, $n\in\N_0$ form an ONS in $L^2(\R_+)$. Let $P(2n+i+\frac{1}{2})$ denote the orthogonal projection onto the one-dimensional subspace spanned by each of these functions.

\begin{The}\label{ALOWF} Let $m=\pm \frac{1}{2}$. Then 
$$T_m\simeq K\oplus K\oplus K\oplus\dots$$
i.e., $T_m$ is Hilbert space isomorphic to the countably infinite orthogonal sum of copies of $K$. More precisely,
$$T_m=\bigoplus_{j\in \N_0+1/2} C_j\;\otimes \,P_m^j $$
where
$\delta C_j \delta^{-1}$ equals $K-\textrm{\SMALL{$\frac{1}{2}$}}\sum_{k=0}^{j-1/2}P(k+\textrm{\SMALL{$\frac{1}{2}$}})$, which is unitarily equivalent to $K$.
\end{The}

{\it Proof.}  The ingredients for the proof are 
\begin{itemize}
\item[(i)] $\big(\mathcal{F}1_B\big)(p-p')=(2/\pi)^{1/2}\int_0^1r^2\operatorname{sinc}(r|p-p'|)\d r$
\item[(ii)] $\operatorname{sinc}(r|p-p'|)=\frac{\pi}{2}\sum_{l=0}^\infty(2l+1)(rsrs')^{-1/2}\operatorname{J}_{l+1/2}(rs) \operatorname{J}_{l+1/2}(rs')\,\operatorname{P}_l\big(\cos\sphericalangle(p,p')
\big)$
\item[(iii)] $(2l+1)\operatorname{P}_l=l\operatorname{P}_{l-1}^{(0,1)}+(l+1)\operatorname{P}_l^{(0,1)}$
\item[(iv)] $D^{(n+1/2)}_{mm}(B^{-1}_wB_{w'})=D^{(1/2)}_{mm}(B^{-1}_wB_{w'})\,\operatorname{P}_n^{(0,1)}\big(\cos\sphericalangle(w,w')\big)$
\item[(v)] $\mathcal{H}_\nu: L^2(\R_+)\to L^2(\R_+)$, $\nu\in[-\frac{1}{2},\infty[$ the unitary involutory Hankel transformation given by $(\mathcal{H}_\nu g)(s)=\int_{\R_+}(ss')^{1/2}\operatorname{J}_\nu(ss')g(s')\d s'$ if $g\in L^2\cap L^1$ and  in general $\mathcal{H}_\nu g=\operatorname{l.i.m.}_{A\to\infty}\mathcal{H}_\nu (1_{[0,A]}g)$.
\end{itemize}
The first (i) is easily checked. For (ii) see Gegenbauer's addition theorem \cite[7.15(30)]{E53} for $\nu=1/2$.
 Item (iii) is \cite[10.8(36)]{E53} for $\alpha=0$, $\beta=1$, since $\operatorname{P}_l=\operatorname{P}_l^{(0,0)}$.
 As to (iv) one recalls that the matrix elements $D^{(j)}_{\sigma\sigma'}(B)$ can be expressed by the Jacobi polynomials $\operatorname{P}^{(\alpha,\beta)}_n$. In particular 
$D^{(j)}_{\sigma\sigma}(B)$ equals $B_{11}^{2|\sigma|}\operatorname{P}^{(0,2|\sigma|)}_{j-|\sigma|}\big(2|B_{11}|^2-1\big)$ or  its complex conjugate depending on whether $\sigma\ge 0$ or $\sigma\le 0$. By this formula one obtains (iv). \\
\hspace*{6mm}
It suffices to evaluate (\ref{CLOUB}) for $\varphi=g\times f $ with integrable $g\in L^2_\rho(\R_+)$, $f\in L^2_\omega(S^2)$. Employing (ii), (iii), and (iv),  (\ref{CLOUB}) becomes 
$\big(T_m\varphi\big)(p)=\int_0^\infty\d s'\int_{S^2}\d \omega(w')\int_0^1\d r \times \sum_{l=0}^\infty\, r^2\, F_l(rs,rs')\,\big(l\,D^{(l-1/2)}_{mm}(B^{-1}_wB_{w'})\,+\,(l+1)\,D^{(l+1/2)}_{mm}(B^{-1}_wB_{w'})\big)\,s'^2\,g(s')\,f(w')$, where
$$ F_l(x,x'):=(xx')^{-1/2}\operatorname{J}_{l+1/2}(x)\operatorname{J}_{l+1/2}(x'), \quad x,x'\in\R_+$$
The integrand is dominated by $\sum_{l=0}^\infty (2l+1)|F_l(rs,rs')|\,1_{[0,1]}(r)\,s'^2g(s')f(w')$. Formula (ii) yields the special case $2/\pi=\sum_{l=0}^\infty (2l+1)F_l(x,x)$,
whence $\sum_{l=0}^\infty (2l+1)|F_l(rs,rs')|\le 2/\pi$ independent of $r,s,s'$. Hence summation and  integrations may be interchanged. Doing first the integral for $\omega$ one obtains $\big(T_m\varphi\big)(p)=\sum_{l=0}^\infty\int_0^1\d r\int_0^\infty\d s' r^2\frac{1}{2}\big(F_l(rs,rs')+F_{l+1}(rs,rs')\big)s'^2g(s')\,(P^{l+1/2}_mf)(w)$. By (iv)  the remaining integrations yield
$T_mg\times f=\sum_{l=0}^\infty \delta^{-1}K_{l+1/2}\delta \,g\otimes P^{l+1/2}_mf$ with
\begin{equation}\label{GOK}
K_{l+1/2}:=\textrm{\SMALL{$\frac{1}{2}$}}\big(\mathcal{H}_{l+1/2}M(1_{[0,1]})\mathcal{H}_{l+1/2}+\mathcal{H}_{l+3/2}M(1_{[0,1]})\mathcal{H}_{l+3/2}\big)
\end{equation} 
Here $M(1_{[0,1]})$ denotes the multiplication operator in $L^2(\R_+)$.\\
\hspace*{6mm}
 $K_{l+1/2}$, $l\in\N\cup\{-1,0\}$ is an integral operator on $L^2(\R_+)$ with kernel 
$$k_{l+1/2}(r,t)=\textrm{\SMALL{$\frac{1}{2}$}}\frac{(rt)^{1/2}}{r-t}\big(\operatorname{J}_{l+3/2}(r)\operatorname{J}_{l+1/2}(t)-\operatorname{J}_{l+1/2}(r)\operatorname{J}_{l+3/2}(t)\big)$$
The latter follows evaluating $(rt)^{1/2}\int_0^1s\operatorname{J}_j(rs)\operatorname{J}_j(st)\d s$ for $j=l+1/2$ and $j=l+3/2$ by \cite[7.14.1(8)]{E53}. Note that $(r-t)k_{l+1/2}(r,t)$ is bounded on $\R_+\times\R_+$ due to \cite[7.2(2)]{E53} for small and \cite[7.13.1(3)]{E53} for large variable $s$ of $s^{1/2}\operatorname{J}_{l+1/2}(s)$. So by \cite[Theorem 6.17.1]{Pu67} the operator 
$K_{l+1/2}$ is absolutely continuous. Now (\ref{FAGKO})(i),(iii) will show $K=K_{-1/2}$ and
\begin{equation}\label{FDPOK}
K_{l+1/2}=K-\textrm{\SMALL{$\frac{1}{2}$}}\sum_{k=0}^{l}P(k+\textrm{\SMALL{$\frac{1}{2}$}})
\end{equation}
which is a a finite-dimensional perturbation of $K$. Hence by \cite[Theorem 6.20.1]{Pu67} all $\delta C_j \delta^{-1}=K_{l+1/2}$ are unitarily equivalent to $K$.\qed

Recall  that $\mathcal{H}_{-1/2}=\mathcal{F}_c$ and  $\mathcal{H}_{1/2} =\mathcal{F}_s$ are the involutory Fourier cosine and  Fourier sine transformation. 

\begin{Lem}\label{FAGKO} 
Let  $n\in\N_0$, $i\in\{0,1\}$. Then
\begin{itemize}
\item[\emph{(i)}] $K=\textrm{\SMALL{$\frac{1}{2}$}}\big(\mathcal{F}_cM(1_{[0,1]})\mathcal{F}_c+\mathcal{F}_sM(1_{[0,1]})\mathcal{F}_s\big)$
\item[\emph{(ii)}] $\mathcal{F}_cM(1_{[0,1]})\mathcal{F}_c=\bigoplus_{n\in\N_0} P(2n+1/2)$,\quad $\mathcal{F}_sM(1_{[0,1]})\mathcal{F}_s=\bigoplus_{n\in\N_0} P(2n+3/2)$
\item[\emph{(iii)}] $\mathcal{H}_{2n+i+1/2}M(1_{[0,1]})\mathcal{H}_{2n+i+1/2}=\mathcal{H}_{1/2-i}M(1_{[0,1]})\mathcal{H}_{1/2-i}-\bigoplus_{k=0}^{n-1+i}P(2k-i+3/2)$
\end{itemize}
\end{Lem}
{\it Proof.} (i) is the special case $K_{-1/2}=K$ of (\ref{GOK}) for $l=-1$ as $k_{-1/2}(r,t)=\frac{1}{\pi}\operatorname{sinc}(r-t)$. As to (ii) specialize the formula (ii) in the proof of (\ref{ALOWF}) to the case of parallel and antiparallel $p,p'$. Note $\operatorname{P}_l(1)=1$, $\operatorname{P}_l(-1)=(-1)^l$. One obtains \begin{itemize}
\item $\frac{1}{\pi}\big(\operatorname{sinc}(x-x')+\operatorname{sinc}(x+x')\big)=\sum_{l=0}^\infty(4l+1)F_{2l}(x,x')$\\
\item $\frac{1}{\pi}\big(\operatorname{sinc}(x-x')-\operatorname{sinc}(x+x')\big)=\sum_{l=0}^\infty(4l+3)F_{2l+1}(x,x')$
\end{itemize}
(for $F_l$ see the proof of (\ref{ALOWF})).
This is (ii) for the  kernels of the relative integral operators.

\hspace*{6mm}
 For the proof of (iii) we first deduce the  formula
 \begin{itemize}
 \item[($\alpha$)] $\frac{\d }{\d z}\sum_{k=n}^{n'}2(\nu+2k)(ab)^{-1}\operatorname{J}_{\nu+2k}(az)\operatorname{J}_{\nu+2k}(bz)=$
 \item[] \begin{flushright}
 $z\big(\operatorname{J}_{\nu+2n-1}(az)\operatorname{J}_{\nu+2n-1}(bz)-\operatorname{J}_{\nu+2n'+1}(az)\operatorname{J}_{\nu+2n'+1}(bz)\big)$
 \end{flushright}
 \end{itemize}
valid for $\nu, z, a,b\in\C$, $ab\ne 0$, $n,n'\in\mathbb{Z}$, $n\le n'$ using the relation
 \begin{itemize}
 \item[($\beta$)] $2\nu\big(z_1\operatorname{J}_{\nu}'(z_1)\operatorname{J}_{\nu}(z_2)+z_2\operatorname{J}_{\nu}(z_1)\operatorname{J}_{\nu}'(z_2)\big)=
 z_1z_2\big(\operatorname{J}_{\nu-1}(z_1)\operatorname{J}_{\nu-1}(z_2)-\operatorname{J}_{\nu+1}(z_1)\operatorname{J}_{\nu+1}(z_2)\big)$
 \end{itemize}
which follows from \cite[7.2.8(54),(55),(56)]{E53} as $2\nu\big(z_1\operatorname{J}_{\nu}'(z_1)\operatorname{J}_{\nu}(z_2)+z_2\operatorname{J}_{\nu}(z_1)\operatorname{J}_{\nu}'(z_2)\big)=2\nu\Big(\big(z_1\operatorname{J}_{\nu-1}(z_1)-\nu\operatorname{J}_{\nu}(z_1)\big)\operatorname{J}_{\nu}(z_2)+\operatorname{J}_{\nu}(z_1)\big(-z_2\operatorname{J}_{\nu+1}(z_2)+\nu\operatorname{J}_{\nu}(z_2)\big)\Big)=2\nu\big(z_1\operatorname{J}_{\nu-1}(z_1)\operatorname{J}_{\nu}(z_2)-z_2\operatorname{J}_{\nu}(z_1)\operatorname{J}_{\nu+1}(z_2)\big)=z_1z_2\Big(\operatorname{J}_{\nu-1}(z_1)\big(\operatorname{J}_{\nu-1}(z_2)+\operatorname{J}_{\nu+1}(z_2)\big)-\big(\operatorname{J}_{\nu-1}(z_1)+\operatorname{J}_{\nu+1}(z_1)\big)\operatorname{J}_{\nu+1}(z_2)\Big)=z_1z_2\big(\operatorname{J}_{\nu-1}(z_1)\operatorname{J}_{\nu-1}(z_2)-\operatorname{J}_{\nu+1}(z_1)\operatorname{J}_{\nu+1}(z_2)\big)$. Now, to the left hand side of ($\alpha$) apply the product rule of differentiation and subsequently use ($\beta$).\\
\hspace*{6mm}
In ($\alpha$) choose $n=0$, $n'=n-1+i$, $\nu=3/2-i$, $a=r$, $b=t$, $z=s$. Integrate ($\alpha$) with respect to $s\in[0,1]$ and note $\operatorname{J}_{\nu}(0)=0$ for $\nu\in\R_+$ according to \cite[7.2.1(2)]{E53}. One obtains the formula (iii) for the  kernels of the relative integral operators.\qed

\subsection{Toeplitz operators isomorphic to $K$}\label{TOIK} Using well-known transformations like that indicated in  \cite[chap.\,9]{BS90} (cf.\,also \cite[sec.\,4]{R65}) we obtain   in (\ref{HSITO}) Toeplitz operators, which are  Hilbert space   isomorphic to $K$, and  in (\ref{SRTO}) their spectral representation. The latter has been  already derived in \cite{R61} for the corresponding discrete Wiener-Hopf operators. In (\ref{GOLBV}) the integral transformation effecting the spectral representation is shown to map the standard basis of the Hardy space
onto the orthonormal basis (\ref{MPPCV}) related to the Meixner-Pollaczek polynomials by a change of variable.

\subsubsection{Toeplitz operators}
 Let $\mathbb{T}\subset \C$ denote the unit circle endowed with the normalized Lebesgue measure. The Hardy space $H^2(\mathbb{T})$ is the closed subspace of $L^2(\mathbb{T})$ with orthonormal basis $e_n,\ n\in\N_0$, where $e_n(z):=z^n$. Let $Q_+:L^2(\mathbb{T})\to H^2(\mathbb{T})$ be the orthogonal projection. The adjoint $Q_+^*$ is the inclusion $H^2(\mathbb{T})\hookrightarrow L^2(\mathbb{T})$. Then for any $\omega\in L^\infty(\mathbb{T})$ consider the multiplication operator 
$M(\omega)$ by $\omega$ in $L^2(\mathbb{T})$ and define the Toeplitz operator 
\begin{equation}\label{TOETS}
T_\omega:=Q_+\,M(\omega)\,Q_+^*
\end{equation}
in $H^2(\mathbb{T})$. Occasionally  we write $T(\omega)$ instead of $T_\omega$.

\subsubsection{Hilbert space isomorphisms transforming Wiener-Hopf operators into Toeplitz operators}
The M\"{o}bius transformations mapping $\R$ to  $\mathbb{T}$ are $\chi(x):=\e  ^{\i  \alpha}\frac{x+\overline{d}}{x+d}$ for $\alpha\in\R$, $d\in\C\setminus\R$, as e.g.\,the Cayley transformation for $\alpha=0$ and $d=\i  $. Each transformation gives rise to a Hilbert space isomorphism $X:L^2(\R)\to L^2(\mathbb{T})$,
\begin{equation}\label{HSIRHT} 
(Xf)(z):=\,(2\pi)^{1/2}\e  ^{\i  \alpha}\,\textrm{\SMALL{$\frac{\sqrt{2|\operatorname{Im} d|}}{-z+\e  ^{\i  \alpha}}$}} f\big(\chi^{-1}(z)\big)
\end{equation}
 Note $\chi^{-1}(z)=\frac{dz-
\e  ^{\i  \alpha}\overline{d}}{-z+\e  ^{\i  \alpha}}$. Its inverse  is
$(X^{-1} h)(x)=(2\pi)^{-1/2}\i  \operatorname{sgn}(\operatorname{Im} d)\,\frac{\sqrt{2|\operatorname{Im} d|}}{x+d}\,h\big(\chi(x)\big)$.    It is easy to verify 
$$XM(\kappa)X^{-1}=M(\kappa\circ\chi^{-1})$$
for $\kappa\in L^\infty(\R)$.

 Recall that $P:=P^*_+P_+=M(1_{[0,\infty[})$ is the orthogonal projection in $L^2(\R)$ with range $L^2(\R_+)$. 
 Similarly $Q:=Q^*_+Q_+$ is the orthogonal projection in $L^2(\mathbb{T})$ with range $H^2(\mathbb{T})$.
  
As to  (\ref{MTRTOWHO}) see \cite[9.1(e),\,(f), 9.2(d), 9.5(e)]{BS90} regarding $U=-\i  X^{-1}$ for $\alpha=\pi$, $d=\i  $.

\begin{The}\label{MTRTOWHO}  Let $\operatorname{Im} d>0$. Then $X^{-1}QX=\mathcal{F}^{-1}P\mathcal{F}$ holds and $X_+: L^2(\R_+)\to H^2(\mathbb{T})$,  $X_+:=Q_+X\mathcal{F}^{-1}P_+^*$
is a Hilbert space isomorphism satisfying $$T(\kappa\circ \chi^{-1})=X_+W_\kappa X_+^{-1}$$
for $\kappa\in L^{\infty}(\R)$.
\end{The}\\
{\it Proof.} Put $f_n:=X^{-1}e_n$ for $n\in\mathbb{Z}$. Then  $X^{-1}QXf_n$ equals $f_n$ if $n\ge 0$ and is $0$ for $n<0$. We are going to show that $\mathcal{F}f_n$ vanishes on $]-\infty,0[$ for  $n\ge 0$ and vanishes on $]0,\infty[$ for $n< 0$, thus implying the first part of the assertion. Note $f_n(x)= \frac{C}{x+d}\chi(x)^n$ with $C$ constant. Hence $\big(\mathcal{F}f_n\big)(x)$ equals the Cauchy integral $C(2\pi)^{-1/2}\int_{-\infty}^\infty \e  ^{-\i  xy} \frac{1}{y+d}\chi(y)^n\d y$, which after a linear change of variable becomes $C'\int_{-\infty}^\infty \e  ^{-\i  \delta xy} \frac{1}{y+i}\big(\frac{y-\i  }{y+\i  }\big)^n\d y$ with $\delta:=\operatorname{Im} d>0$ and $C':=(2\pi)^{-1/2}\e  ^{\i  (n\alpha+x\gamma)}$, $\gamma:=
\operatorname{Re} d$. Note that for $n\ge 0$ there are no residues in the upper half-plane and  that for $n<0$ there are none in the lower-half plane. For $x<0$ the integral can be evaluated closing the path along the real line from $-A$ to $A$ by the half-circle in the upper half-plane and performing the limit $A\to\infty$. Similarly for $x>0$ one closes the path  by the half-circle in the lower half-plane. So the result holds by the residue theorem.\\
\hspace*{6mm}
Hence one has $QX\mathcal{F}^{-1}=X\mathcal{F}^{-1}P$, whence  $X_+^*X_+=I_{L^2(\R_+)}$,  $X_+X_+^*=I_{H^2(\mathbb{T}}$  and the last part of the assertion.\qed

Recall that $l_q(x):= (1-q)^{-1}\exp\big(-\frac{1}{2}\frac{1+q}{1-q}\, x\big)=\sum_{n=0}^\infty \operatorname{l}_n(x)q^n$, $|q|<1$ is  the generating function of the Laguerre functions $l_n$, $n\ge 0$, which form an orthonormal basis in $L^2(\R_+)$. Similarly $e_q(z):=\frac{1}{1-qz}$  generates  the standard basis $(e_n)_{n\ge 0}$ of $H^2(\mathbb{T})$.

As to (\ref{ILFUIXP}) cf.\,\cite[sec.\,4,\,iv)]{R65} regarding $X_+$ for $\alpha=0$, $d=\frac{\i  }{2}$.

\begin{Lem}\label{ILFUIXP} Consider $X$ in \emph{(\ref{HSIRHT})} for $\operatorname{Im} d>0$. If $\alpha=0$, $d=\frac{\i  }{2}$
one has  $X_+l_n=e_n$, $n\ge 0$. Generally $$X_+l_n(z)=\i  \e  ^{\i  \alpha}\sqrt{2\operatorname{Im} d}\,\big((d-\textrm{\SMALL{$\frac{\i  }{2}$}})z-\e  ^{\i  \alpha}(\overline{d}-\textrm{\SMALL{$\frac{\i  }{2}$}})\big)^{-1}e_n\big(\mu^{-1}(z)\big)$$ holds, where the M\"obius transformation $\mu(z):=\e  ^{\i  \alpha}\frac{(\overline{d}-\textrm{\SMALL{$\frac{\i  }{2}$}})z-\overline{d}-\textrm{\SMALL{$\frac{\i  }{2}$}}}{(d-\textrm{\SMALL{$\frac{\i  }{2}$}})z-d-\textrm{\SMALL{$\frac{\i  }{2}$}}}$ maps $\mathbb{T}$ onto itself.
\end{Lem}\\
{\it Proof.} One easily computes $\mathcal{F}^{-1}P_+^*l_q(x)=(2\pi)^{-1/2}\frac{\i  }{x+\frac{\i  }{2}}\,e_q\Big(\frac{x-\frac{\i  }{2}}{x+\frac{\i  }{2}}\Big)=X_0^{-1}e_q(x)$, where the subscript $0$ denotes the case $\alpha=0, d=\frac{\i  }{2}$. 
Therefore $X_+l_n=Q_+X\mathcal{F}^{-1}P_+^*l_n=Q_+XX_0^{-1}e_n$, $n\ge 0$, whence the result.\qed

\subsubsection{Toeplitz operators isomorphic to $K$ and their spectral representations}\label{TOITK} Let $A$  denote the  arc $\{\e  ^{\i  \varphi}:\alpha+\beta\le\varphi\le\alpha+2\pi-\beta\}$ 
for  $\alpha\in [0,2\pi[$, $\beta\in]0,\pi[$. 
One easily verifies $1_A=1_{[-1,1]}\circ\chi^{-1}$ for $\chi(x)=\e  ^{\i  \alpha}\frac{x-\i  \tan(\beta/2)}{x+\i  \tan(\beta/2)}$, \,$\chi^{-1}(z)= \i  \tan(\beta/2)\frac{z+\e  ^{\i  \alpha}}{-z+\e  ^{\i  \alpha}}$,
whence by (\ref{MTRTOWHO}) 
\begin{equation}\label{HSITO} 
 T(1_A)=X_+\,K\,X_+^{-1}
\end{equation}
 
Recall the Hilbert space isomorphism $V=V_+$ in (\ref{HSISFK}) providing the spectral representation  of $K$  (\ref{KHSIJ}). 
Then obviously the Hilbert space isomorphism $Y:=X_+V: L^2(0,1)\to H^2(\mathbb{T})$ yields the spectral representation 
\begin{equation}\label{SRTO} 
 T(1_A)=Y\,M(\operatorname{id}_{[0,1]})\,Y^{-1}
\end{equation}

\hspace*{6mm} For $n\in\N_0$ let $h_n:]0,1[\to \R$
 \begin{equation}\label{MPPCV}
h_n(s):=\e  ^{\i  n\alpha}\Big(\textrm{\SMALL{$\frac{\sin\beta}{\pi}$}}\Big)^{1/2}\big(\textrm{\SMALL{$\frac{1}{s}$}}-1\big)^{\beta/(2\pi)}(1-s)^{-1/2}\operatorname{P}_n^{(1/2)}\Big(\textrm{\SMALL{$\frac{1}{2\pi}$}}\ln\big(\textrm{\SMALL{$\frac{1}{s}$}}-1\big);\beta\Big)
\end{equation}
where $\operatorname{P}_n^{(1/2)}(x;\beta)$ is the $n$th Meixner-Pollaczek polynomial of order $1/2$.  The generating function of the  polynomials is $p_q(x):=(1-q\e  ^{\i  \beta})^{-1/2 + \i  x}(1-q\e  ^{-\i  \beta})^{-1/2 - \i  x}$. Orthogonality   and completeness of the polynomials $\operatorname{P}_n^{(1/2)}$ on $\R$ carry over  to the functions $h_n$ by the change of variable $x=\frac{1}{2\pi}\ln\big(\frac{1}{s}-1\big)$ so that the latter form an orthonormal basis of $L^2(0,1)$.

As to (\ref{GOLBV}) cf.\;the result in \cite[Example 3]{R61}.

\begin{The}\label{GOLBV}\quad  $Y^{-1}e_n=h_n$,   $n\in\N_0$.
\end{The}\\
{\it Proof.} Note $Y^{-1}=V^{-1}X_+^{-1}$. Use the abbreviation $y:=\frac{1}{2\pi}\ln\big(\textrm{\SMALL{$\frac{1}{s}$}}-1\big)$ and put $t(u):=\ln\big(\frac{1-u}{1+u}\big)$. Then by (\ref{KSFK}), (\ref{EKSFK}), and due to (\ref{MPQS})(i) one has $\sqrt{2}\,\pi s\sqrt{1-s}\,(Y^{-1}h)(s)=(2\pi)^{-1/2}\operatorname{l.i.m.}\int_{-\infty}^{\uparrow\infty}\Big(\int_{-1}^1\e  ^{-\i  xu}(1-u^2)^{-1/2}\e  ^{\i  yt(u)}\d u\Big)(\mathcal{F}X^{-1}h)(x)\d x$.  Let the transform $X^{-1}h$ of $h\in H^2(\mathbb{T})$ have compact support. Then the integrations can be interchanged. Use $\mathcal{F}^2f(u)=f(-u)$. It follows  for all $h\in H^2(\mathbb{T})$
\begin{equation*}(Y^{-1}h)(s)=\big(\sqrt{ 2}\,\pi \,s\sqrt{1-s}\big)^{-1}  \operatorname{l.i.m.} \textrm{\SMALL{$\int_{\downarrow-1}^{\uparrow 1}$}}(1-u^2)^{-1/2}\e  ^{-\i  yt(u)}(X^{-1}h)(u)\d u \tag{$\star$}
\end{equation*}
Now we evaluate ($\star$) for $h=e_q$.  We treat the case $\alpha=0$. For the general case simply substitute  $q$ with $q\e  ^{\i  \alpha}$. Put $\delta:=\tan\frac{\beta}{2}>0$ and $C(s):=\i  \sqrt{2\delta}\,\big( 2\,\pi^{3/2} \,s\sqrt{1-s}\big)^{-1}$. Perform the change of variable $t=t(u)$. One obtains $Y^{-1}e_q(s)=C(s)\int_{-\infty}^\infty \e  ^{-\i  yt}\frac{\e  ^{t/2}}{c_q-\overline{c}_q\e  ^t}\d t$ for $c_q:=1-q+\i  \delta(1+q)$ and $q\in]-1,1[$.  Note $\overline{c}_q/c_q=\e  ^{\i  \gamma}$ for some $\gamma\in]-\pi,0[$. Hence the residue theorem yields for all $y\in\R$
$$\int_{-\infty}^\infty \e  ^{-\i  yt}\frac{\e  ^{t/2}}{1-\e  ^{\i  \gamma}\e  ^t}\d t=-2\pi\i  \frac{\e  ^{-\gamma(y+\i  /2)}}{1+\e  	^{2\pi y}}$$
Simple computations show  
$c_q=\sec\beta/2\,\e  ^{\i  \beta/2}\big(1-q\e  ^{-\i  \beta}\big)$, $\e  ^{\i  \gamma}=\e  ^{-\i  \beta}
\frac{1-q\e  ^{\i  \beta}}{1-q\e  ^{-\i  \beta}}$, whence 
$\e  ^{-\gamma(y+\i  /2)}=\Big(\frac{1-q\e  ^{\i  \beta}}{1-q\e  ^{-\i  \beta}}\Big)^{-1/2+\i  y}\e  ^{\i  \beta/2}\e  ^{\beta y}$. Also note that  $|\operatorname{arg}(1+w)|\le\frac{\pi}{2}$ if $|w|<1$. Therefore putting together the terms one obtains $Y^{-1}e_q(s)=\big(\frac{\sin\beta}{\pi}\big)^{1/2}\big(\frac{1}{s}-1\big)^{\beta/(2\pi)}(1-s)^{-1/2}p_q\big(\frac{1}{2\pi}\ln(\frac{1}{s}-1)\big)$, whence the result.\qed\\

 Performing in ($\star$) the change of variable $z=\chi(u)$  one gets the integral transformation 
 \begin{equation}
 (Y^{-1}h)(s)=\i  (2\tan\beta/2)^{1/2} (\pi s)^{-1/2}\int_A\frac{Q\big(1-2s,\chi^{-1}(z)\big)}{-z+\e  ^{\i  \alpha}}h(z)\d z
 \end{equation}
with $Q$ from (\ref{IKFHT}) for $h\in H^2(\mathbb{T})$ with compact support within the open arc. Because of (\ref{SRTO}) the kernel of this transformation yields  normalized spectral functions of $T(1_A)$ for every spectral value $s\in]0,1[$.\\

\hspace*{6mm}
Consider the special case $\alpha=0$ and $\tan\frac{\beta}{2}=\frac{1}{2}$. Then (\ref{GOLBV}) yields by (\ref{ILFUIXP})
\begin{equation}\label{SGOLBV}
V^{-1}\operatorname{l}_n=h_n 
\end{equation}
and, more precisely, 
$ \int_0^\infty q_+(s,x)\operatorname{l}_n(x)\d x=h_n(s) \quad\forall \;s\in]0,1[,\; n\in\N_0$.

\subsection{Finite Hilbert transformation}\label{MEFHT} The following consideration show
a close relationship  between  $K$ and the finite Hilbert transformation on $[-1,1]$.

\subsubsection{Hilbert transformation} The spectral representation of the  finite Hilbert transformation on $L^2(a,b)$, $-\infty<a<b<\infty$
\begin{equation}\label{FHT}
(H_{[a,b]}\,g)(x):=\textrm{\SMALL{$\frac{1}{\i  \pi}$}}\int_{a}^b\frac{g(y)}{y-x}\d y
\end{equation} 
(in the sense of the principal value at $x$) is given in  \cite{KP59}.  See also \cite[Theorem 6.20.1]{Pu67} for another approach.   
  In \cite[Lemma 2.1]{KP59}   spectral functions of $H_{[a,b]}$ are found by solving a related boundary value problem. Then, by the Hilbert space isomorphism $U_{a,b}$  associated to the former obtained in \cite[Theorem 3.1]{KP59},   the representation
\begin{equation}\label{SRFHT}
H_{[a,b]}=\,U_{a,b}\,M(\operatorname{id}_{[-1,1]})\,U_{a,b}^{-1}
\end{equation}
is shown in \cite[Theorem 3.2]{KP59}). Here $M(\operatorname{id}_{[-1,1]})$ denotes the multiplication operator on $L^2(-1,1)$ by the identity.\\
\hspace*{6mm}
Recall  the well-known  representation of the  Hilbert transformation  $H$ ($a=-\infty$, $b=\infty$)
\begin{equation}\label{SPHT}  
H=\mathcal{F}^{-1}\,  M(\operatorname{sgn})\,\mathcal{F}= -\mathcal{F}\,  M(\operatorname{sgn})\,\mathcal{F}^{-1}
\end{equation}
where $M(\operatorname{sgn})$ denotes the multiplication operator on $L^2(\R)$ by the signum function. See e.g.\,\cite[Theorems 91,\,95]{T48} or \cite[Lemma 1.35]{D79}, and for more details \cite[Teorema 1.1.1]{T73}. Hence
\begin{equation}\label{FHTFR}
H_{[a,b]}=\,P_{[a,b]}\,\mathcal{F}^{-1}\,M(\operatorname{sgn})\,\mathcal{F}\,P_{[a,b]}^*
\end{equation}
where $P_{[a,b]}:L^2(\R)\to L^2(a,b)$ denotes the projection $\big(P_{[a,b]}f\big)(u):=f(u)$. Its adjoint $P_{[a,b]}^*$ is the injection $\big(P_{[a,b]}^*k\big)(x)=k(x)$ if $a\le x\le b$ and $0$ otherwise. Analogous formulae hold for the semi-finite transformations  $H_{]-\infty,b]}$ and $H_{[a,\infty[}$. \\
 \hspace*{6mm}
In the following $a=-1$, $b=1$ is assumed. This is no restriction as $H_{[a,b]}$  is obtained by an obvious linear change of variable. Even the semi-finite transformations are Hilbert space isomorphic to the finite one using a fractional linear change of variable as shown in \cite[sec.\,4]{KP59}. Indeed, apply  the Hilbert space isomorphism $\Gamma:L^2(a,b)\to L^2(-1,1)$, $(\Gamma g)(t):=\sqrt{|\gamma'(t)|}\,g(\gamma(t))$, $\gamma:]-1,1[\to ]a,b[$ with $\gamma(t):= \frac{b-a}{2}t+  \frac{b+a}{2}$ if $a,b$ are finite, and $\gamma(t):=b+\frac{t-1}{t+1}$ if $a=-\infty$, $\gamma(t):=a+\frac{1+t}{1-t}$ if $b=\infty$. 

\subsubsection{Spectral analysis of $H_{[-1,1]}$} 
 
  We are going to derive the spectral analysis of $H_{[-1,1]}$ from  the foregoing analysis of $K$, comprising  the spectral representation 
 (\ref{SRFHT}), which is  shown in (\ref{PSRFHT}). \\
 \hspace*{6mm}
Introduce $A:=P_+\mathcal{F}P_{[-1,1]}^*$, $A^*=P_{[-1,1]}\mathcal{F}^{-1}P_+^*$. Using (\ref{FPK}), (\ref{FHTFR}), and $M(1_{[0,\infty[})=\frac{1}{2}\big(1+M(\operatorname{sgn})\big)$
one easily checks \begin{itemize}
\item[(a)] $A\,A^*=K$
\item[(b)] $A^*A=\frac{1}{2}(I+H_{[-1,1]})$
\end{itemize}

The relation (a), (b) is not accidental and as shown in \cite{C20} it relates isomorphicly   
 self-adjoint Wiener-Hopf operators to singular integral operators with kernels of Cauchy type. \\
 \hspace*{6mm}
The following computations are valid in a  distributional sense.  Spectral functions $Q_s$ for $\frac{1}{2}(I+H_{[-1,1]})$  are given by $A^*q_s =P_{[-1,1]}\mathcal{F}^{-1}P_+^*q_s$ since 
 $(A^*A)\,A^*q_s=A^*(AA^*)q_s=sA^*q_s$. By  (\ref{MPQS})(i),\,(iii),\,(i) one gets explicitly
 \begin{equation}\label{SFIPHP} 
Q_s=\frac{n(s)}{\sqrt{s}}\; Q^+_s|_{[-1,1]}, \quad s\in]0,1[
\end{equation}
Here $\frac{n(s)}{\sqrt{s}}=\i  (\sqrt{2}\,\pi s)^{-1}$ is the unique (up to a phase) normalization constant such that $Q(s,x):=Q_s(x)$ is  the kernel for a Hilbert space isomorphism  $W:L^2(0,1)\to L^2(-1,1)$ satisfying
\begin{equation}\label{HSIDIPHP} 
 Wh=\operatorname{l.i.m.}\int_ {\downarrow 0}^{\uparrow 1} Q(s,\cdot)h(s)\d s
\end{equation} 
 The additional factor $\frac{1}{\sqrt{s}}$ regarding the normalization  of $Q_s$  is suggested heuristically by $\langle A^*q_s, A^*q_s\rangle=\langle q_s,AA^*q_s\rangle=s\langle q_s,q_s\rangle$. Indeed, 
 (\ref{HSIDIPHP}) follow immediately from 
 \begin{equation}\label{ADHSIDIPHP}
W= \Gamma\,\mathcal{F}^{-1}\Gamma^{-1}\Lambda
\end{equation}
 for the Hilbert space isomorphisms
$\Lambda:L^2(0,1)\to L^2(-1,1)$, $\Lambda h\,(t):=\frac{1}{\sqrt{2}} h\big(\frac{1+t}{2}\big)$
due to the change of variable $\lambda(t):=(1+t)/2$, and $\Gamma: L^2(\R)\to L^2(-1,1)$,  $(\Gamma f)(x):=\sqrt{|\gamma'(x)|}f(\gamma(x))$ due to the change of variable $\gamma:]-1,1[\to\R$, $\gamma(x):=\sqrt{\frac{1}{2\pi}} \ln\big(\frac{1+x}{1-x}\big)$, $\gamma'(x)=\big(\frac{2}{\pi}\big)^{1/2}(1-x^2)^{-1}$. In order to verify  (\ref{ADHSIDIPHP}) one easily checks that the unitary transformation 
\begin{equation}
U:=\Gamma\,\mathcal{F}^{-1}\Gamma^{-1}
\end{equation}
satisfies
\begin{equation}
Uk=\operatorname{l.i.m.}_{0<\alpha,\alpha'\to 0}\int_{-1+\alpha}^{1-\alpha'}Q'(t,\cdot)k(t)\d t
\end{equation}
for $Q'(t,u):=\frac{1}{\sqrt{2}}Q\big(\frac{1+t}{2}\big)$. By (\ref{SFIPHP}),\,(\ref{LQS}) one has explicitly
\begin{equation}\label{IKFHT}
Q'(t,u)=\textrm{\SMALL{$\frac{1}{\pi}$}}(1-t^2)^{-1/2}(1-u^2)^{-1/2}\exp\textrm{\SMALL{$\left(\frac{\i  }{2\pi}\ln\Big(\frac{1-t}{1+t}\Big) \ln \Big(\frac{1-u}{1+u}\Big)\right)$}}
\end{equation}

Now we give another proof of  \cite[Theorem 3.2]{KP59}).

\begin{Cor}\label{PSRFHT} The  spectral representation 
$H_{[-1,1]}=U\,M(\operatorname{id}_{[-1,1]})\,U^{-1}$ holds.
\end{Cor}\\
{\it Proof.}  It suffices to show $P_{[-1,1]}\mathcal{F}^{-1}\,M(1_{[0,\infty[})\,\mathcal{F}P_{[-1,1]}^*Wh=W\,M(\operatorname{id}_{[0,1]})h$ for  $h$ with compact support in $]0,1[$. For $A>1$, 
$\int_{-A}^A\e  ^{-\i  ux}\big(\int_{0}^11_{]-1,1[}(u)Q(s,u)h(s)\d s\big)\d u=\int_{0}^1\big(\int_{-1}^1\e  ^{-\i  ux}Q(s,u)\d u\big)h(s)\d s$ for every $x\in\R$ by Fubini's theorem since $|Q(s,u)h(s)|$ is integrable by (\ref{LQS}). By (\ref{MPQS})(ii) the double integral becomes $\i  \pi^{-1/2}\int_{0}^1q_s^+(x)h(s)\d s$. The latter is zero if $x<0$ by (\ref{MPQS})(i), and due to (\ref{MPQS})(iii) equals $\int_{-\infty}^\infty\e  ^{-\i  ux}\Big(\int_{0}^1Q(s,u)sh(s) \d s\Big) \d u$. This implies the assertion.\qed\\

\subsection{Covariance under real fractional linear transformations} One starts from the following simple observation. Recall that $H$ denotes the Hilbert transformation.

\begin{Lem}\label{CVUT} Let $F:L^2(\R)\to L^2(\R)$ be unitary.  Put $F':=P_+\mathcal{F}F\mathcal{F}^{-1}P_+^*$.  Suppose
\begin{itemize}
\item[\emph{($\alpha$)}] $FH=HF$
\item[\emph{($\beta$)}]
For every $\kappa\in L^\infty(\R)$ there is $ \kappa'\in L^\infty(\R)$ such that $FM(\kappa)F^{-1}=M(\kappa')$
\end{itemize}
Then $F'$ is unitary and $F'W_\kappa F'^{-1}=W(\kappa')$
holds.
\end{Lem}\\
{\it Proof.} By (\ref{SPHT}) and ($\alpha$),  $M(1_{[0,\infty[})$ commutes with $\mathcal{F}F\mathcal{F}^{-1}$ reducing it to  $F'$. Hence $F'P_+=P_+\mathcal{F}F\mathcal{F}^{-1}$, whence the result by ($\beta$).\qed\\

 With $A=\big(\textrm{\SMALL{$ \begin{array}{cc} a & b\\ c & d \end{array}$}}\big) \in SL(2,\R)$
  associate the real fractional linear transformation $A\cdot x:=\frac{ax+b}{cx+d}$ for almost all  $x\in\R$. 
 Actually this is a group action as $A\cdot (A'\cdot x)= (AA' )\cdot x$ and $I_2\cdot x=x$. Note that $A'\cdot x=A\cdot x$ for almost all $x$  if and only if $A'\in\{A,-A\}$.
  
\begin{The}\label{CVFLTWH}  For $A\in SL(2,\R)$, $f\in L^2(\R)$ let $\big(F_Af\big)(x):=\frac{1}{-cx+a}f(A^{-1}\cdot x)$. Then $F$ is a unitary representation of $SL(2,\R)$ in $L^2(\R)$ commuting with $H$. Let $F'$ denote the subrepresentation of $\mathcal{F}F\mathcal{F}^{-1}$ according to \emph{(\ref{CVUT})}. One has the covariance
$$W(A\cdot \kappa)\,=\,F'_A\,W_\kappa\, F_A'^{-1}$$
where $\big(A\cdot \kappa\big)(x):=\kappa(A^{-1}\cdot x)$ for $\kappa\in L^\infty(\R)$. 
 \end{The}\\
 {\it Proof.} The result follows by elementary computations. Verify $F_AH=HF_A$ and $M(A\cdot \kappa)=F_AM(\kappa)F_A^{-1}$, and apply (\ref{CVUT}).\qed\\
 
 In particular   $W(1_J)$ for $J$ a bounded or semi-bounded interval is unitarily equivalent to $K$ by some $F'_A$. For instance $A:=\frac{1}{\sqrt{2}}\big(\textrm{\SMALL{$ \begin{array}{cc} 1-\gamma & 1+\gamma\\ -1 & 1 \end{array}$}}\big)$ provides for the isomorphism of $K$ to $W(1_{[\gamma,\infty[})$, which has the distribution kernel $k(x)=\frac{1}{2}\delta(x)+ \frac{1}{2\pi \i  }\frac{\exp(-\i  \gamma x)}{-x}$, whence $W(1_{[0,\infty[})=\frac{1}{2}(I+H_{[0,\infty[})$ if $\gamma=0$ (cf.\,(\ref{FHTFR})). More generally, $W(a1_{]-\infty,\gamma[}+b1_{[\gamma,\infty[})$
 has the kernel $k(x)=\frac{b+a}{2}\delta(x)+ \frac{b-a}{2\pi \i  }\frac{\exp(-\i  \gamma x)}{-x}$.
 \\

For (\ref{TLWHHT})($\beta$) see the remark on  (0.3) in \cite{D79}.

\begin{Lem} \label{TLWHHT} One has
\begin{itemize} 
\item[\emph{($\alpha$)}] $H_{[0,\infty[}=W_{-\operatorname{sgn}}$, where 
$W_{-\operatorname{sgn}}$ has the kernel $-\frac{1}{x}$ 
\item[\emph{($\beta$)}] $H_{[0,1]}=\Gamma \,W_{-\tanh}\Gamma^{-1}$, where $W_{-\tanh}$ has the kernel $\frac{1}{2\i  \sinh(\pi x/2)}$  
\end{itemize} 
Here  $\Gamma:L^2(\R_+)\to L^2(0,1)$, $(\Gamma g)(t):=(\pi t)^{-1/2}g\big(-\frac{1}{\pi}\ln(t)\big)$ is the Hilbert space isomorphism due to the change of variable $\gamma(t):=-\frac{1}{\pi}\ln(t)$.
\end{Lem}

{\it Proof.}  For ($\alpha$) see (\ref{FHTFR}). As to ($\beta$), check $(\Gamma^{-1}k)(x)=\sqrt{\pi}\e  ^{-\pi x/2}k(\e  ^{-\pi x})$ and  that $\Gamma^{-1}H_{[0,1]}\Gamma$ equals the Wiener-Hopf operator with the indicated kernel. Then its symbol $-\tanh(x)$ is the principal value of $\int \frac{\exp(\i  xy)}{2\i  \sinh(\pi y/2)}\d y$ at $y=0$, which is evaluated calculating the residues.\qed

We use these results in an obvious way  to relate explicitly $K$ to several other operators. Note that $2K-I=W_{\sigma}$ with $\sigma(x)=1$ if $|x|\le1$ and $=-1$ otherwise.

\begin{Cor}\label{HTWHRK} The finite and semi-finite Hilbert transformations  as well the Wiener-Hopf operators $W_\kappa$ for $\kappa(x)=-\tanh(\frac{ax+b}{cx+d})$ with $ad-bc=1$ are Hilbert space isomorphic to $2K-I$. More precisely
$$H_{[a,b]}=\Lambda\,(2K-I)\,\Lambda^{-1}, \quad W_\kappa=\Lambda'\,(2K-I)\,\Lambda'^{-1}$$
Here $\Lambda= \Gamma\, F'_B$, where  $B:=\frac{1}{\sqrt{2}}\big(\textrm{\SMALL{$ \begin{array}{cc} 1 & -1\\ 1 & 1 \end{array}$}}\big)$ and $\Gamma:L^2(0,\infty)\to L^2(a,b)$ due to the change of variable $\gamma(t):=\frac{t-a}{-t+b}$ if $a,b$ are finite, $\gamma(t):=t-a$ if $b=\infty$, $\gamma(t):=\frac{1}{-t+b}$ if $a=-\infty$. Further 
 $\Lambda':=F_A'^{-1}\Gamma^{-1}\Lambda$ with $A=\big(\textrm{\SMALL{$ \begin{array}{cc} a & b\\ c & d \end{array}$}}\big)$, $\Gamma$ from \emph{(\ref{TLWHHT})($\beta$)}, and $\Lambda$ relating $H_{[0,1]}$ to $2K-I$.
 \end{Cor}\\
 {\it Proof.} Put $F:=F_B$. Check $B\cdot (21_{]-1,1[}-1)=-\operatorname{sgn}$, whence $F'(2K+I)F'^{-1}=H_{[0,\infty[}$ by (\ref{TLWHHT})($\alpha$),\,(\ref{CVFLTWH}). Further check $\Gamma H_{[0,\infty[}\Gamma^{-1}=H_{[a,b]}$. Hence the first part of the assertion holds. The remainder follows from this result for $H_{[a,b]}=H_{[0,1]}$ and by (\ref{TLWHHT})($\beta$),  (\ref{CVFLTWH}).\qed\\
 
 In particular (\ref{HTWHRK}) yields
\begin{equation} \label{CHSIHK}
H_{[-1,1]}\,=\,\Lambda (2K-I) \Lambda^{-1}\; \textrm{ with }  \Lambda=P_{[-1,1]}F_B\mathcal{F}F_B\mathcal{F}^{-1}P_+^*
 \end{equation}

\end{document}